\documentstyle[grcalc]{elsart}

\def\myauthor{Ph\`ung H{{\accent"5E o}\kern-.28em\raise.2ex\hbox{\char'22}\kern-.20em}
 H{a\kern-.370em\raise.16ex\hbox{\char'47}\kern.1em}i}

\def\iii{^{-1}}

\def\lora{\longrightarrow}
\def\ot{\otimes}

\def\si{\sigma}
\newcommand{\bbas}{\begin{eqnarray*}}
\newcommand{\eeas}{\end{eqnarray*}}
\newcommand{\bbar}{\begin{array}}
\newcommand{\eear}{\end{array}}
\newcommand{\bbs}{\begin{displaymath}}
\newcommand{\ees}{\end{displaymath}}
\newcommand{\bb}{\begin{equation}}
\def\ee{\end{equation}}
\def\eea{\end{eqnarray}}

\def\bba{\begin{eqnarray}}

\newtheorem{edl}[thm]{Theorem}

\def\coend{\mbox{\rm coend}}
\def\Mor{\mbox{\rm Mor}}

\def\C{\mbox{$\cal C$}}

\def\V{{\cal V}}

\def\D{{\cal D}}
\def\E{{\cal E}}

\def\eee{\mbox{\rule{.75ex}{1.5ex}}}

\newcommand{\proof}{{\it Proof.\ }}
\newcommand{\va}{\varepsilon}
\def\nml{\normalsize}

\newcommand{\Coalg}{\mbox{\rm Coalg}}

\def\rref#1{(\ref{#1})}
\newcommand{\M}[1]{\mbox{${\cal M}^#1$}}

\font\Fraktur=eufm10 scaled\magstep1          % for display- and textstyle
   \newcommand{\fraktur}[1]{\mbox{\Fraktur #1}}  %
   \font\Fraktu=eufm7 scaled\magstep1            % for scriptstyle
   \newcommand{\fraktu}[1]{\mbox{\Fraktu #1}}    %
   \font\Frakt=eufm5 scaled\magstep1             % for scriptscriptstyle
  \newcommand{\frakt}[1]{\mbox{\Frakt #1}}      %
   \def\fr#1{\mathchoice{\fraktur {#1}}            % displaystyle
                        {\fraktur {#1}}            % textstyle
                        {\fraktu {#1}}             % scriptstyle
                        {\frakt {#1}}  }           % scriptscriptstyle

\newcommand{\Cc}{\fr C}
\newcommand{\Bb}{\fr B}
\newcommand{\Aa}{\fr A}
\newcommand{\Zz}{\fr Z}

\font\tes=eusm10 scaled \magstep1

\newcommand{\test}[1]{\mbox{\tes #1}}

\newcommand{\Ccc}{{\test C}}
\newcommand{\Bbb}{{\test B}}
\newcommand{\Aaa}{{\test A}}

\def\db{{\mathchoice{\mbox{\rm db}}
                    {\mbox{\rm db}}
                    {\mbox{\scriptsize\rm db}}
                    {\mbox{\tiny\rm db}} }}

\def\ev{{\mathchoice{\mbox{\rm ev}}
                    {\mbox{\rm ev}}
                    {\mbox{\scriptsize\rm ev}}
                    {\mbox{\tiny\rm ev}} }}

\def\id{{\mathchoice{\mbox{\rm id}}
                    {\mbox{\rm id}}
                    {\mbox{\scriptsize\rm id}}
                    {\mbox{\tiny\rm id}} }}

\def\op{{\mathchoice{\mbox{\rm op}}
                    {\mbox{\rm op}}
                    {\mbox{\scriptsize\rm op}}
                    {\mbox{\tiny\rm op}} }}

\begin{document}
\begin{frontmatter}
\title{ Central Bialgebras in Braided Categories and Coquasitriangular Structures}
\author{\myauthor}

\address{International Center for Theoretical Physics,\\ P.O. Box 586, 34100 Trieste, Italy \\ email: phung@ictp.trieste.it }

\begin{abstract} Central bialgebras in a braided category $\C$ are algebras in the center of the category of coalgebras in $\C$. On these bialgebras another product can be defined, which plays the role of the opposite product. Hence, coquasitriangular structures on central bialgebras can be defined. We prove some properties of the antipode on  coquasitriangular central Hopf algebras and give a characterization of central bialgebras.\end{abstract}

\end{frontmatter}

\footnotetext[1]{ JOURNAL OF PURE AND APPLIED ALGEBRA to appear}

\section*{Introduction}
In his talk at the ICM 1986 V. Drinfel'd introduced the notion of quantum groups \cite{drinfeld86}. By definition, quantum groups are Hopf algebras which own a distinguished element in their tensor square, called $R$-matrix, which controls the cocommutativity of the Hopf algebras. Also of interest is the dual notion, called cotriangle, or in more general case dual quasitriangular or coquasitriangular Hopf algebras, first appeared in the works of Lyubashenko, Majid \cite{lyu,x}. These are Hopf algebras $H$, which are quasicocomutative in the sense that the opposite product coincides with the original product up to conjugation by a linear functional $r:H\ot H\lora k$, where $k$ is the ground field. They can also be characterized \cite{x,y} as such that the category $\M{H}$ of right $H$-comodules is a braided category, i.e., a monoidal category which possesses a natural isomorphism, called braiding, which twists the tensor product of every two objects. Braided categories themselves were introduced in category theory by Joyal and Street \cite{js1}.

Recently many authors are interested in braided groups, which are defined to be Hopf algebras in  braided categories. The notions of algebras and coalgebras can still be nicely generalized to braided categories. However trouble appears when one studies bialgebras in braided categories. There are two main obstructions in the theory of bialgebras in braided categories. First, the (co) opposite (co) product defined in the standard way does not satisfy its axioms. Second,  the tensor product of two bialgebra, in a braided category, does not have a bialgebra structure. Fortunately, there is no trouble with the notion of the antipode, since it definition does not involve the braiding. Lacking the (co) opposite (co) product one can not define the notion of (co) commutativity of bialgebras in braided categories.

To solve  the first  obstruction mentioned above, Majid \cite{majid1} suggested considering bialgebras with a second product, which is also a coalgebra morphism. This second product is called weak opposite product in the sense that it is opposite (to the original one)  with respect to a class of comodules. Having two products one can then define the notion of commutativity, coquasitriagularity (with respect to a class of comodules). Examples of this construction are bialgebras reconstructed from monoidal functors, on which there is a natural choice for the second product.

The obstructions in the theory of bialgebras in braided categories mentioned above can be explained in the following way. Observe, that in a monoidal category, one can define the notion of algebras and coalgebras but not the notion of (co) commutative (co) algebras. Neither can one define a (co) algebra structure on the tensor product of two (co) algebras, in other words, the category of (co) algebras in a monoidal category is not monoidal. One solves this problem by introducing a braiding. Thus, the category of (co) algebras in a braided category is monoidal. On the other hand, one can define the notion of the opposite (co) product for (co) algebras in a braided category. A bialgebra in a braided category $\C$ can be considered as an algebra in the category of coalgebras in $\C$ or a coalgebra in the category of algebras in $\C$.  Since the category of (co) algebras in a braided category is generally not braided, the category of bialgebras in braided categories in not monoidal. That explains why we cannot define a bialgebra structure on the tensor product of two bialgebras as well as the notion of the opposite (co) product on a bialgebra. 

Of course we cannot ``make'' a monoidal category into a braided category. But there is a standard way of obtaining a braided category from a monoidal category, that is taking the center of this category \cite{js2,y}. Roughly speaking, an object of the center of a monoidal category is an object of the last category that has additional properties. However, it is to emphases that an object of a monoidal category may appear in many ways to be an object of the center of this category, thus, the center of a monoidal category is far from a subcategory of this category. Subject of the present paper is the category of algebra in the center of the monoidal category of coalgebras in a braided category. This category is therefore monoidal. Its objects, which we call central bialgebras, can be seen as bialgebras in the initial braided category with some extra properties. It is then not surprising, that for central bialgebras one can define in a natural way the second product, that plays the role of the opposite product -- one uses the new braiding. To summarize, instead of trying to solve the obstructions for every bialgebra, we shrink the class of bialgebras being considered, so that our theory on this class will be more consistent.

In Section \ref{sec-central} we define  central bialgebras in a braided category as algebras in the center of the category of coalgebras in this category. For this class of bialgebras the opposite product can be defined in a very natural way. In Section \ref{sec-cqt} we define coquasitriangular structures on central bialgebras and study the antipode of a coquasitriangular central Hopf algebra. We describe a method of obtaining central bialgebras in Section \ref{sec-exam}.
\subsection*{Preliminary}\label{sec-pre}
Let $\C$ be a category. A monoidal structure in $\C$ consists of a bifunctor $\ot:\C\times \C\lora \C,$ an object $I\in \C$ and natural   isomorphisms $\lambda:I\ot C\lora  C,\mu: C\ot I\lora C, C\in\C$, $\Phi_{A,B,C}:(A\ot B)\ot C\lora A\ot(B\ot C)$ one for every triple $A,B,C\in \C$, satisfying certain coherence-conditions \cite{maclane}. A monoidal category is called strict if the morphism $\lambda, \mu, \Phi$ are identity-morphisms.  Later on we will assume all monoidal categories considered here to be strict.

Let $(\C,\ot)$ be a monoidal category. A braiding in $\C$ is a natural isomorphism
\bbas \tau_{A,B}:A\ot B\lora B\ot A,\eeas
 satisfying the following equations:
\bb\label{braidings}\bbar{c} \tau_{A\ot B,C}=(\tau_{A,C}\ot B)(A\ot \tau_{B,C}),\\
\tau_{A,B\ot C}=(B\ot \tau_{A,C})(\tau_{A,B}\ot C),\eear\ee
here (and later on) we use the same letter to denote an object and the identity morphism on it.

We will frequently use the graphical calculus, our main references are \cite{js1,majid2}. In particular a morphism $f:A\lora B$ and the braiding $\tau, \tau^{-1}$ will be depicted as follows\cite{majid2}.

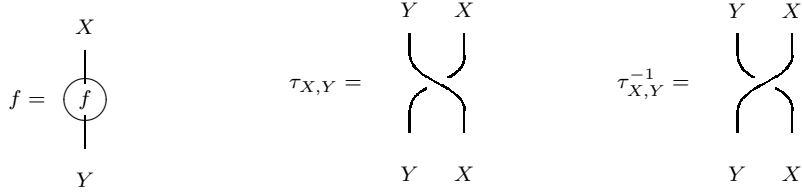
\cbgr{130}{20}\put(-10,10){$f=$}
\put(4,0){\obju{Y}}\put(4,20){\objo{X}}\put(4,2){\mor{f}{18}}
\put(50,0){\bgr{.4ex}{70}{22}
\put(12,0){\obju{Y}}\put(22,0){\obju{X}}\put(12,22){\objo{Y}}\put(22,22){\objo{X}}
\put(12,4){\idgr{4}}\put(22,4){\idgr{4}}\put(12,18){\idgr{4}}\put(22,18){\idgr{4}}\put(12,8){\braid}\put(-10,12){$\tau_{X,Y}=$}
\put(72,0){\obju{Y}}\put(82,0){\obju{X}}\put(72,22){\objo{Y}}\put(82,22){\objo{X}}
\put(72,4){\idgr{4}}\put(82,4){\idgr{4}}\put(72,18){\idgr{4}}\put(82,18){\idgr{4}}\put(72,8){\ibraid}\put(50,12){$\tau_{X,Y}^{-1}=$}\egr}
\cegr

In a monoidal category one can define the notion of algebras, coalgebras. In a braided category one can define the (co) algebra structure on  the tensor product of (co) algebras hence the notion of bialgebras and Hopf algebras. Through out this paper, $m$ and $\eta$ (resp. $\Delta$ and $\va$) will denote the product and unit of an algebra (resp. the coproduct and counit of a coalgebra). In graphical notation the (co) product and (co) unit are depicted as follows:

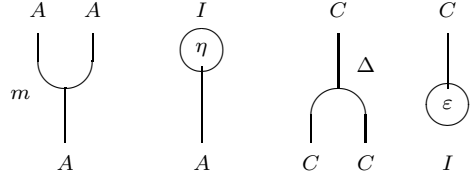
\cbgr{80}{20}
 \put(0,10){\mult}\put(5,0){\idgr{10}}\put(0,15){\idgr{5}}\put(10.05,15){\idgr{5}}\put(0,20){\objo{A}}
\put(10,20){\objo{A}}\put(5,0){\obju{A}}\put(-5,8){$m$}\put(30,0){\umor{\eta}{20}}\put(30,20){\objo{I}}\put(30,0){\obju{A}}

\put(50,0){\idgr{5}}\put(60,0){\idgr{5}}
\put(55,10){\idgr{10}}
\put(50,5){\comult}\put(60,10){\objo{\Delta}}
\put(75,4){\coumor{\varepsilon}{16}}\put(75,0){\obju{I}}
\put(50,0){\obju{C}}
\put(60,0){\obju{C}}
\put(55,20){\objo{C}}
\put(75,20){\objo{C}}\cegr

The antipode on a Hopf algebra is denoted by $S$.

For a coalgebra $C$, $\Mor(C,I)$ is a monoid with the product, denoted by $*$, defined as follows:
\bbas -*-=(-\ot -)\Delta:C\lora I.\eeas
The unit in $\Mor(C,I)$ is $\va_C$. We call this product convolution- or $*$-product.

Let $(\C,\ot)$ be a monoidal category. The center of $\C$ \cite{js2,y}, denoted by $\Zz(\C)$, is a category, the objects of which are pairs $(A,\si_A(-))$ consisting of an object $A$ of $\C$ and a natural isomorphism in $\C$ $\si_A(N): A\ot N\lora N\ot A$, satisfying
\bbas \si_A(I)=\id_A,\\
 \si_A(M\ot N)= (M\ot \si_A(N))(\si_A(M)\ot N).\eeas
 Morphisms in $\Zz(\C)$ are those in $\C$ which commute with $\si$. $\Zz(\C)$ is a braided category. For  objects $(A,\si_A(-))$ and $(B,\si_B(-))$ in $\Zz(\C)$ their tensor product is defined to be $(A\ot B, (\si_A\ot B)(A\ot \si_B))$ and the braiding is given by $\si_A(B)$, we will denote it also by $\si_{A,B}$ when no confusion may arise.

\section{Central bialgebras}\label{sec-central}
Let $(\C,\tau)$ be a braided category,  algebras (resp. coalgebras) in $\C$ form a monoidal category denoted by $\Aa(\C)$ (resp. $\Cc(\C)$). Morphisms in $\Aa(\C)$ (resp. $\Cc(\C)$ ) are morphism of algebras (resp. coalgebras) in $\C$. 

For algebras $A,B$ in $\C$, the algebra structure on  their tensor product is given by means of the braiding in $\C$: $m_{A\ot B}:=(m_A\ot m_B)(A\ot \tau_{B,A}\ot B)$. Analogously, for coalgebras $A$ and $B$ we have $\Delta_{A\ot B}:=(A\ot\tau_{B,A}\ot B)(\Delta_A\ot \Delta_B).$ The categories $\Aa(\C)$ and $\Cc(\C)$ are are not braided unless the braiding in $\C$ is symmetric. 

The category $\Bb(\C)$ of bialgebras in $\C$ coincides with the categories $\Aa(\Cc(\C))$ and $\Cc(\Aa(\C))$. $\Bb(\C)$  is, in general, even not monoidal. This is the main obstruction in the theory of bialgebras in braided category. In order to get more structure in the category we attack as follows. Since $\Cc(\C)$ is monoidal, $\Zz(\Cc(\C))$ is a braided category. Hence $\Aa(\Zz(\Cc(\C)))$ is a monoidal category.  Analogously, one can define the monoidal category $\Cc(\Zz(\Aa(\C)))$. We call objects of $\Aa(\Zz(\Cc(\C)))$ central bialgebras and denote this category by $\Zz\Bb(\C)$. Objects of $\Cc(\Zz(\Aa(\C)))$ are called cocentral bialgebras. This paper is focused on central bialgebras. To get similar results for cocentral bialgebras  rotate  all the pictures 180 degrees around the horizontal axe.

Let me recall the axioms of $\Zz\Bb(\C)$. Objects of $\Zz\Bb(\C)$ are pairs $(B,\si_B)$ consisting of a bialgebra $B$ and a natural isomorphism $\si_B(-)$ in $\Cc(\C)$, such that the multiplication and the unit of $B$ commute with $\si_B$. Thus $\si_B$ satisfies the following axioms:

\begin{description}
\item[Z1] $\si_{B,C}$ is a coalgebra morphism, for all coalgebra $C$ in $\C$,
\bb \label{eqz1}\Delta_{C\ot B} \si_{B,C}= (\si_{B,C}\ot \si_{B,C})\Delta_{B\ot C},\ee
\item[Z2] $\si_{B,C}$ is natural in $C$:
\bb\label{eqz2}\si_{B,C}(B\ot f)=(f\ot B) \si_B(D),\ee
 for all coalgebra morphism $f:C\lora D$,
\item[Z3] 
\label{eqz3}There exists $ \si^{-1}_B(C):C\ot B\lora B\ot C,$ for all coalgebra $ C$,
\item[Z4] The unit $\eta$ of $B$ commutes with $\si_{B,C}$:
\bb\label{eqz4}\si_{B,C} (\eta\ot C)=C\ot \eta,\ee
\item[Z5] The product $m$ of $B$ commutes with $\si_{B,C}$:
\bb\label{eqz5} \si_{B,C} m=(C\ot m) \si_{B\ot B,C}.\ee 
  \end{description}

The morphism $\si_{B,B}$ is called the innertwist of $B$ and will be denoted briefly by $\si$ if now confusion may arise.

In graphical notations \rref{eqz1}--\rref{eqz5} are depicted as follows, (where  {\scriptsize\begin{minipage}{6mm}\bgr{.3ex}{5}{5}\put(0,2){\brmor}\egr\end{minipage}} denotes $\si$),

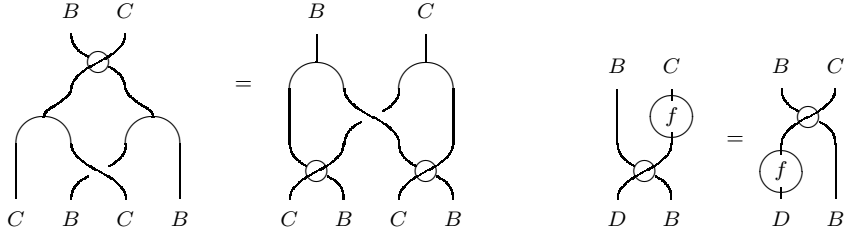
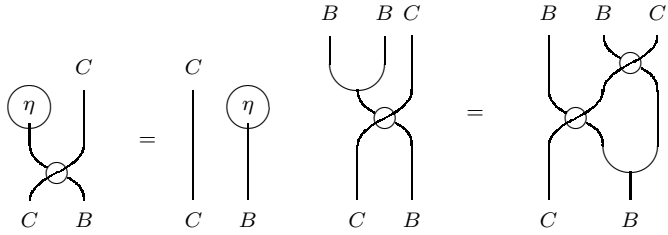
\cbgr{150}{30}\put(0,0){\bgr{.4ex}{20}{35}\put(0,0){\idgr{10}}\put(10,0){\braid}\put(30,0){\idgr{10}}\put(0,10){\comult}\put(20,10){\comult}\put(5,15){\twist{5}{5}}\put(25,15){\twist{-5}{5}}\put(10,20){\simor}\put(0,0){\obju{C}}\put(20,0){\obju{C}}\put(10,0){\obju{B}}\put(30,0){\obju{B}}\put(10,30){\objo{B}}\put(20,30){\objo{C}}\put(40,20){=}\egr}
\put(50,0){\bgr{.4ex}{20}{35}\put(0,0){\simor}\put(20,0){\simor}\put(10,10){\braid}\put(0,10){\idgr{10}}\put(30,10){\idgr{10}}\put(0,20){\comult}\put(20,20){\comult}\put(5,25){\idgr{5}}\put(25,25){\idgr{5}}\put(0,0){\obju{C}}\put(20,0){\obju{C}}\put(10,0){\obju{B}}\put(30,0){\obju{B}}\put(5,30){\objo{B}}\put(25,30){\objo{C}}\egr}

\put(110,0){\bgr{.4ex}{20}{35}\put(0,0){\simor}\put(0,10){\idgr{10}}\put(10,10){\mor{f}{10}}\put(20,10){=}\put(0,0){\obju{D}}\put(10,0){\obju{B}}\put(0,20){\objo{B}}\put(10,20){\objo{C}}\egr}
\put(140,0){\bgr{.4ex}{20}{35}\put(0,0){\mor{f}{10}}\put(10,0){\idgr{10}}\put(0,10){\simor}\put(0,0){\obju{D}}\put(10,0){\obju{B}}\put(0,20){\objo{B}}\put(10,20){\objo{C}}\egr}\cegr

\cbgr{140}{40}
\put(0,0){\bgr{.4ex}{20}{35}\put(0,0){\simor}\put(10,10){\idgr{10}}\put(0,10){\umor{\eta}{10}}\put(20,10){=}\put(0,0){\obju{C}}\put(10,0){\obju{B}}\put(10,20){\objo{C}}\egr}
\put(30,0){\bgr{.4ex}{20}{35}\put(10,0){\umor{\eta}{20}}\put(0,0){\idgr{20}}\put(0,0){\obju{C}}\put(10,0){\obju{B}}\put(0,20){\objo{C}}\egr}

\put(55,0){\bgr{.4ex}{20}{35}\put(5,0){\idgr{10}}\put(15,0){\idgr{10}}\put(5,10){\simor}\put(15,20){\idgr{10}}\put(0,20){\mult}\put(0,25){\idgr{5}}\put(10,25){\idgr{5}}\put(5,0){\obju{C}}\put(15,0){\obju{B}}\put(0,30){\objo{B}}\put(10,30){\objo{B}}\put(15,30){\objo{C}}\put(25,15){=}\egr}
\put(95,0){\bgr{.4ex}{20}{35} \put(0,0){\idgr{10}}\put(15,0){\idgr{5}}\put(0,10){\simor}\put(0,20){\idgr{10}}\put(10,5){\mult}\put(10,20){\simor}\put(20,10){\idgr{10}}\put(0,0){\obju{C}}\put(15,0){\obju{B}}\put(0,30){\objo{B}}\put(10,30){\objo{B}}\put(20,30){\objo{C}}\egr}\cegr

We denote $\si^{-1}$ by {\scriptsize\begin{minipage}{6mm}\bgr{.3ex}{5}{5}\put(0,2){\ibrmor}\egr\end{minipage}}.

For any coalgebra $C$ the counit $\va:C\lora I$ is a morphism of coalgebra. Hence, applying to both sides of \rref{eqz1} the morphism $C\ot B\ot \va\ot B$ and $\va\ot B\ot C\ot B$ respectively, we get:
\bb\label{sigl}\begin{array}{ll}  (C\ot \Delta_B)\si_{B,C}&
= (\si_{B,C}\ot B)(B\ot \tau_{B,C})(\Delta_B\ot C)\\ & =(\tau^{-1}_{C,B}\ot B)(B\ot \si_{B,C})(\Delta_B\ot C).\end{array}\label{eqz11}\ee
Note that these equations together with \rref{eqz2} imply \rref{eqz1}. 

\rref{eqz1} can be reformulated in term of $\si^{-1}$ which derives the following equations:
\bb\label{isigl}\begin{array}{ll}  (\Delta_B\ot C)\si^{-1}_{B,C}& = (B\ot\si^{-1}_{B,C})(\tau_{C,B}\ot B)(C\ot \Delta_B)\\ & = (B\ot \tau^{-1}_{B,C})(\si_{B,C}^{-1}\ot B)(C\ot \Delta_B).\end{array}\label{eqz12}\ee
In graphical notations  \rref{eqz11} and \rref{eqz12} have the following forms:

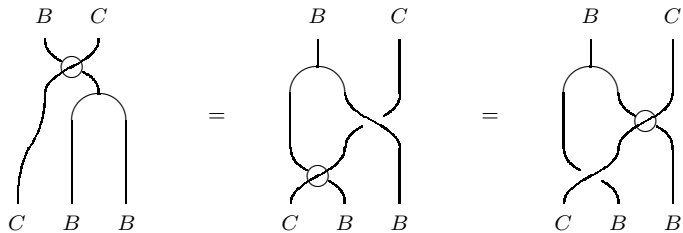
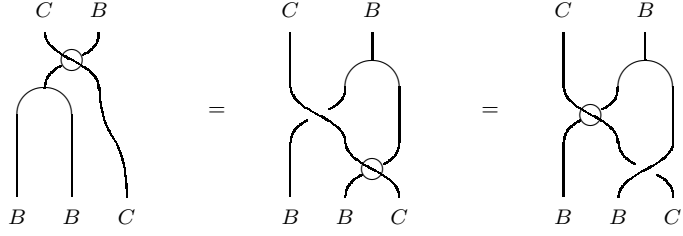
\cbgr{120}{30}
\put(0,0){\bgr{.4ex}{100}{30}\put(0,0){\twist{5}{20}}\put(5,20){\sibrmor}\put(10,0){\idgr{15}}\put(20,0){\idgr{15}}\put(10,15){\comult}\put(0,0){\obju{C}}\put(10,0){\obju{B}}\put(20,0){\obju{B}}\put(05,30){\objo{B}}\put(15,30){\objo{C}}\put(35,15){=}\egr}
\put(50,0){\bgr{.4ex}{100}{30}\put(0,0){\sibrmor}\put(0,10){\idgr{10}}\put(0,20){\comult}\put(10,10){\braid}\put(5,25){\idgr{5}}\put(20,20){\idgr{10}}\put(20,0){\idgr{10}}\put(0,0){\obju{C}}\put(10,0){\obju{B}}\put(20,0){\obju{B}}\put(5,30){\objo{B}}\put(20,30){\objo{C}}\put(35,15){=}\egr}
\put(100,0){\bgr{.4ex}{100}{30}\put(0,0){\ibraid}\put(0,10){\idgr{10}}\put(0,20){\comult}\put(10,10){\sibrmor}\put(5,25){\idgr{5}}\put(20,20){\idgr{10}}\put(20,0){\idgr{10}}\put(0,0){\obju{C}}\put(10,0){\obju{B}}\put(20,0){\obju{B}}\put(5,30){\objo{B}}\put(20,30){\objo{C}}\egr}\cegr

\cbgr{120}{35}
\put(0,0){\bgr{.4ex}{100}{30}\put(0,0){\idgr{15}}\put(10,0){\idgr{15}}\put(20,0){\twist{-5}{20}}\put(0,15){\comult}\put(5,20){\isibrmor}\put(0,0){\obju{B}}\put(10,0){\obju{B}}\put(20,0){\obju{C}}\put(5,30){\objo{C}}\put(15,30){\objo{B}}\put(35,15){=}\egr}
\put(50,0){\bgr{.4ex}{100}{30}\put(0,0){\idgr{10}}\put(0,20){\idgr{10}}\put(0,10){\braid}\put(10,0){\isibrmor}\put(20,10){\idgr{10}}\put(10,20){\comult}\put(15,25){\idgr{5}}\put(0,0){\obju{B}}\put(10,0){\obju{B}}\put(20,0){\obju{C}}\put(0,30){\objo{C}}\put(15,30){\objo{B}}\put(35,15){=}\egr}
\put(100,0){\bgr{.4ex}{100}{30}\put(0,0){\idgr{10}}\put(0,20){\idgr{10}}\put(0,10){\isibrmor}\put(10,0){\ibraid}\put(20,10){\idgr{10}}\put(10,20){\comult}\put(15,25){\idgr{5}}\put(0,0){\obju{B}}\put(10,0){\obju{B}}\put(20,0){\obju{C}}\put(0,30){\objo{C}}\put(15,30){\objo{B}}\egr}\cegr
These two equations will be used very frequently in the next section.
We also remark that, for a bialgebra $C$, the unit $\eta_C$ is a coalgebra morphism. Hence, $\si_{B,C}(B\ot \eta_C)=\eta_C\ot B$.

If $B$ is a Hopf algebra then the antipode, denoted by $S$, is an anti-(co)algebra morphism \cite{majid1}:
\bbas \Delta S= \tau_{B,B}(S\ot S),\\ Sm=m\tau_{B,B}(S\ot S).\eeas

Let  $(C,\si)$ be a bialgebra, consider $C$ as coalgebra, then $C^{\op}$, defined in the standard way by $\Delta^{\op}:=\tau\Delta$, is a coalgebra, too. Hence, $S_C$ is a coalgebra morphism from $C$ to $C^{\op}$, therefore commutes with $\si_{B,C}$:
\bb \si_{B,C}(B\ot S)=(S\ot B)\si_{B,C}.\ee

A central bialgebra $(B,\si)$  with $B$ being a Hopf algebra  is called central Hopf algebra. The alternative equation, in  which $S$ and $B$ change their places in the tensor product, is far from being true. In fact, we have the following lemma, which is due to P. Schauenburg.
\begin{lem} Let $B$ be a central Hopf algebra, then the antipode in $B$ satisfies
\bb\label{sigmaantipode} \si_{B,C}(S\ot C)=(C\ot S)\tau_{B,C}\si^{-1}_{B,C}\tau^{-1}_{C,B}.\ee\end{lem}
\proof Let us consider the morphism

$(C\ot m( S\ot B))
(\tau_{B,C}\si^{-1}_{B,C}\tau^{-1}_{C,B}\ot m)(B\ot\si_{B\ot B,C}) 
(B^{\ot 2}\ot S\ot B) (\Delta^{(2)}_B\ot C)$ \\ $ :B\ot C\lora C\ot B.$

There are two ways of reducing it, which give on one hand $(C\ot S)(\tau_{B,C}\si^{-1}_{B,C}\tau^{-1}_{C,B})$ and on the other hand $\si_{B,B}(S\ot C)$. Let me use the graphical calculus for showing this (the morphism above is depicted by the middle picture).

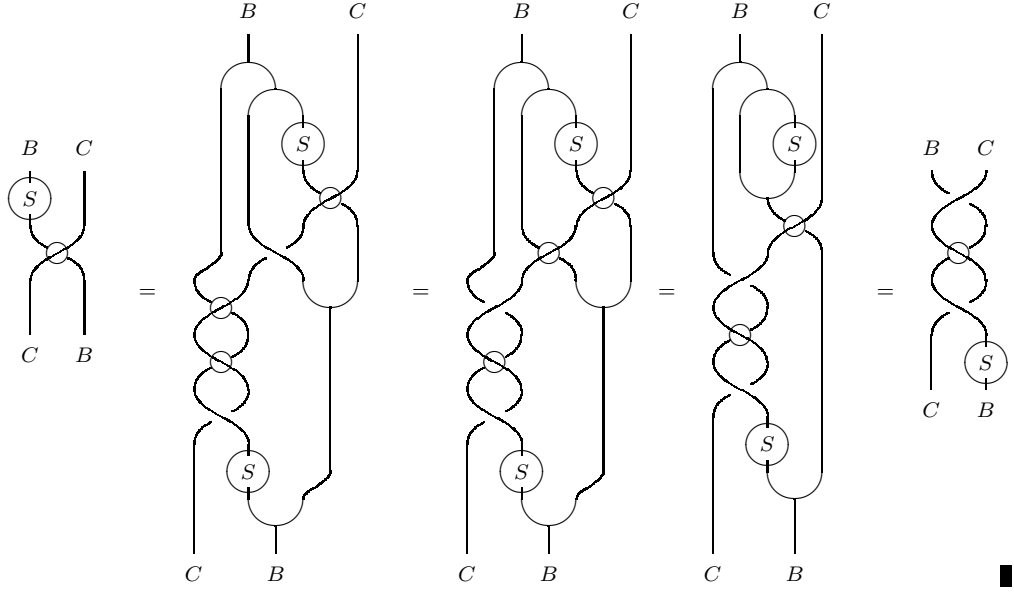
\cbgr{175}{95}\eeee
\put(0,0){\bgr{.4ex}{20}{35}\put(0,40){\idgr{10}}\put(10,40){\idgr{10}}\put(10,60){\idgr{10}}\put(0,60){\mor{S}{10}}\put(0,50){\simor}\put(0,40){\obju{C}}\put(10,40){\obju{B}}\put(0,70){\objo{B}}\put(10,70){\objo{C}}\put(20,47){=}\egr}

\put(30,0){\bgr{.4ex}{20}{35}\put(0,0){\idgr{20}}\put(15,0){\idgr{5}}\put(10,5){\mult}\put(20,10){\twist{5}{5}}\put(10,10){\mor{S}{10}}\put(0,20){\braid}\put(25,15){\idgr{30}}\put(0,30){\isimor}\put(0,40){\simor}\put(0,50){\twist{5}{5}}\put(5,55){\idgr{30}}\put(10,50){\braid}\put(20,45){\mult}\put(30,50){\idgr{10}}\put(10,60){\idgr{20}}\put(20,60){\simor}\put(20,70){\mor{S}{10}}\put(30,70){\idgr{25}}\put(10,80){\comult}\put(5,85){\comult}\put(10,90){\idgr{5}}\put(30,50){\idgr{10}}\put(0,0){\obju{C}}\put(15,0){\obju{B}}\put(10,95){\objo{B}}\put(30,95){\objo{C}}\put(40,47){=}
\egr}

\put(80,0){\bgr{.4ex}{20}{35}\put(0,0){\idgr{20}}\put(15,0){\idgr{5}}\put(10,5){\mult}\put(20,10){\twist{5}{5}}\put(10,10){\mor{S}{10}}\put(0,20){\braid}\put(25,15){\idgr{30}}\put(0,30){\isimor}\put(0,40){\ibraid}\put(0,50){\twist{5}{5}}\put(5,55){\idgr{30}}\put(10,50){\simor}\put(20,45){\mult}\put(30,50){\idgr{10}}\put(10,60){\idgr{20}}\put(20,60){\simor}\put(20,70){\mor{S}{10}}\put(30,70){\idgr{25}}\put(10,80){\comult}\put(5,85){\comult}\put(10,90){\idgr{5}}\put(30,50){\idgr{10}}\put(0,0){\obju{C}}\put(15,0){\obju{B}}\put(10,95){\objo{B}}\put(30,95){\objo{C}}\egr}

\put(125,0){\bgr{.4ex}{20}{35}\put(0,0){\idgr{25}}\put(15,0){\idgr{10}}\put(10,10){\mult}\put(10,15){\mor{S}{10}}\put(20,15){\idgr{40}}\put(0,25){\braid}\put(0,35){\isimor}\put(0,45){\ibraid}\put(10,55){\simor}\put(0,55){\idgr{30}}\put(5,65){\mult}\put(20,65){\idgr{30}}\put(5,70){\idgr{10}}\put(15,70){\mor{S}{10}}\put(5,80){\comult}\put(0,85){\comult}\put(5,90){\idgr{5}}\put(-10,47){=}\put(30,47){=}\put(0,0){\obju{C}}\put(15,0){\obju{B}}\put(5,95){\objo{B}}\put(20,95){\objo{C}}\egr}

\put(165,0){\bgr{.4ex}{20}{35}\put(0,30){\idgr{10}}\put(10,30){\mor{S}{10}}\put(0,40){\braid}\put(0,50){\isimor}\put(0,60){\ibraid}\put(0,30){\obju{C}}\put(10,30){\obju{B}}\put(0,70){\objo{B}}\put(10,70){\objo{C}}\egr}
\cegr

The category $\Zz\Bb(\C)=\Aa(\Zz(\Cc(\C)))$ is monoidal, that means the tensor product of central bialgebras is again a central bialgebra. The coproduct on the tensor product of two central bialgebras is defined by means of the braiding in $\C$, while the product is defined by means of the braiding  $\si$ in $\Zz(\Cc(\C))$: for central bialgebras $B$ and $C$, 
\bbas \Delta_{B\ot C}:=(B\ot\tau_{B,C}\ot C)(\Delta_B\ot \Delta_C),\\
m_{B\ot C}:=(m_B\ot m_C)(B\ot \si_{B,C}\ot C).\eeas
If moreover $B$ and $C$ are Hopf algebras then their tensor product is  a Hopf algebra, too, the antipode is given by \cite{neuchl}:
\bbs S_{B\ot C}=\si_{B,C}(S_C\ot S_B)\tau^{-1}_{B,C}.\ees
If one of the bialgebras (Hopf algebras) is not central their tensor product  still has  a bialgebra (Hopf algebra) structure defined in the same way, but it is no more central \cite{neuchl}.

We now come  to the notion of the opposite product. Let $(B,\si)$ be a central bialgebra. Considering $(B,\si)$ as an algebra in $\Zz(\Cc(\C))$, since $\Zz(\Cc(\C))$ is braided one can define a new  algebra structure by means of the braiding in $\Zz(\Cc(\C))$. Thus  $B^{\op}$ with the product defined by  $m^{\op}:=m\si_{B,B}$ is a central bialgebra.
 The theorem below will deal with the case of central Hopf algebras.
\begin{edl} Let $(B,\si)$ be a central Hopf algebra in $\C$. If the antipode on $B$ is invertible then $(B^{\op},\si)$ is a central Hopf algebra in $\C$. The antipode on $B^{\op}$ is given by:

\bbs \overline{S}=(\va\ot S^{-1})\si^{-1}\tau^{-1}\Delta:B\lora B.\ees
\end{edl}
\proof The following equations shows that $\overline{S}$ obeys the axioms for an antipode:

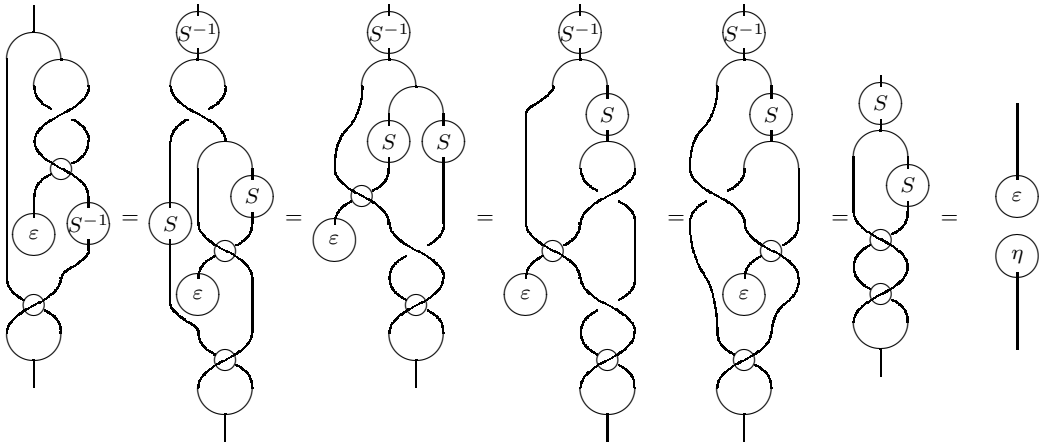
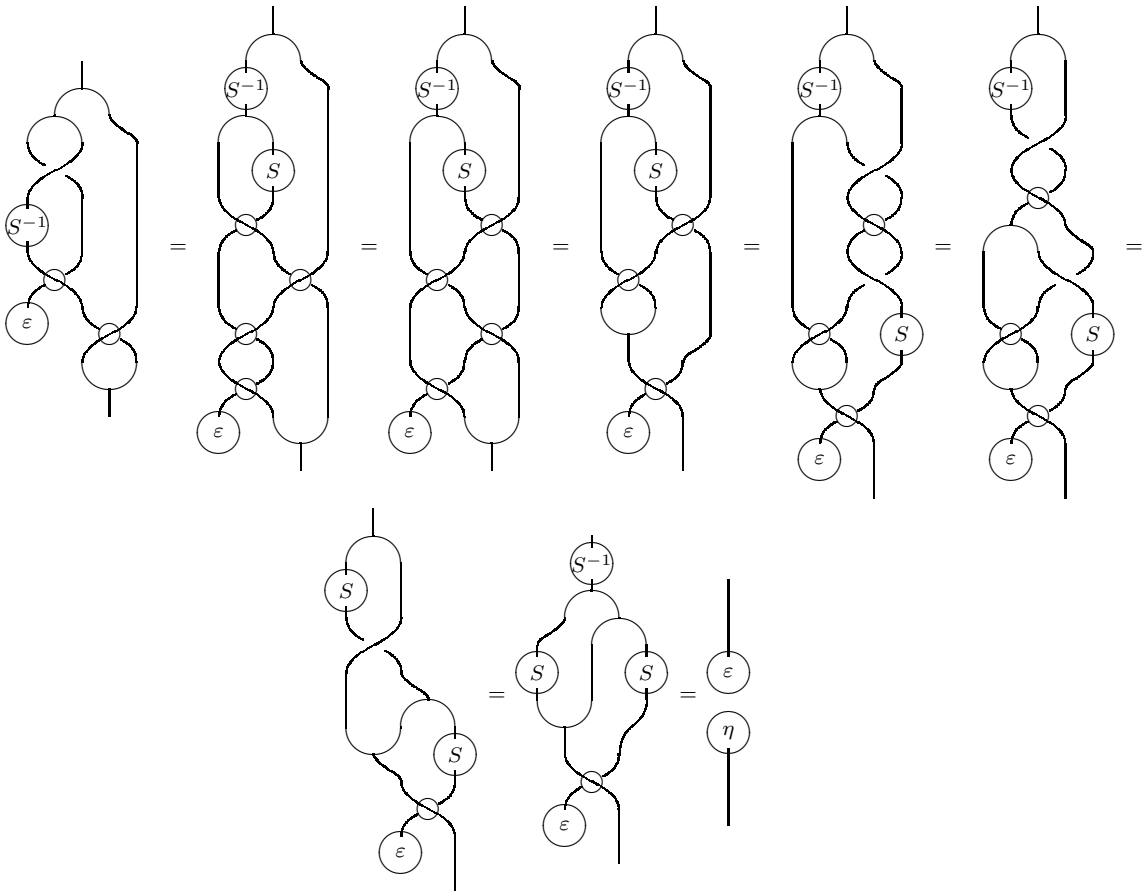
\cbgr{185}{80}
\put(0,10){\bgr{.4ex}{15}{80}\put(05,0){\ig{5}}\put(0,5){\mult}\put(0,10){\sibrmor}\put(0,20){\ig{40}}\put(0,60){\comult}\put(5,65){\ig{5}}\put(10,20){\twist{5}{5}}\put(5,25){\coumor{\va}{10}}\put(5,35){\isibrmor}\put(5,45){\ibraid}\put(5,55){\comult}\put(15,25){\mor{S\iii}{10}}\put(21,30){=}\egr}%%%%%%%%
\put(30,0){\bgr{.4ex}{15}{80}\put(10,0){\ig{5}}\put(5,5){\mult}\put(5,10){\sibrmor}\put(5,20){\twist{-5}{5}}\put(0,25){\mor{S}{30}}\put(0,55){\braid}\put(0,65){\comult}\put(5,70){\mor{S\iii}{10}}\put(5,24){\coumor{\va}{6}}\put(5,30){\isibrmor}\put(5,40){\ig{10}}\put(5,50){\comult}\put(15,40){\mor{S}{10}}\put(15,20){\ig{10}}\put(21,40){=}\egr}%%%%%%

\put(60,10){\bgr{.4ex}{15}{80}\put(15,0){\ig{5}}\put(10,5){\mult}\put(10,10){\sibrmor}\put(10,20){\braid}\put(0,24){\coumor{\va}{6}}\put(0,30){\isibrmor}\put(20,30){\ig{10}}\put(0,40){\twist{5}{15}}\put(5,55){\comult}\put(10,50){\comult}\put(10,60){\mor{S\iii}{10}}\put(10,40){\mor{S}{10}}\put(20,40){\mor{S}{10}}\put(26,30){=}\egr}
\put(95,0){\bgr{.4ex}{20}{80}\put(15,0){\ig{5}}\put(10,5){\mult}\put(10,10){\sibrmor}\put(10,20){\braid}\put(0,24){\coumor{\va}{6}}\put(0,30){\isibrmor}\put(20,30){\ig{10}}\put(0,40){\ig{20}}\put(10,40){\ibraid}\put(10,50){\comult}\put(0,60){\twist{5}{5}}\put(5,65){\comult}\put(15,55){\mor{S}{10}}\put(10,70){\mor{S\iii}{10}}\put(26,40){=}\egr}
\put(125,0){\bgr{.4ex}{20}{80}\put(10,0){\ig{5}}\put(5,5){\mult}\put(5,10){\sibrmor}\put(5,20){\twist{-5}{20}}\put(15,20){\twist{5}{10}}\put(10,24){\coumor{\va}{6}}\put(10,30){\isibrmor}\put(0,40){\braid}\put(10,50){\comult}\put(20,40){\ig{10}}\put(0,50){\twist{5}{15}}\put(5,65){\comult}\put(10,70){\mor{S\iii}{10}}\put(15,55){\mor{S}{10}}\put(26,40){=}\egr}
\put(155,12){\bgr{.4ex}{10}{80}\put(5,0){\ig{5}}\put(0,5){\mult}\put(0,10){\sibrmor}\put(0,20){\isibrmor}\put(0,30){\ig{10}}\put(10,30){\mor{S}{10}}\put(0,40){\comult}\put(5,45){\mor{S}{10}}\put(16,28){=}\egr}
\put(185,17){\bgr{.4ex}{10}{80}\put(0,0){\umor{\eta}{20}}\put(0,25){\coumor{\va}{20}}\egr}\cegr
\noindent and 

\cbgr{200}{75}
\put(0,10){\bgr{.4ex}{10}{85}\put(15,0){\ig{5}}\put(10,5){\mult}\put(10,10){\sibrmor}\put(0,14){\coumor{\va}{6}}\put(0,20){\isibrmor}\put(20,20){\ig{30}}\put(0,30){\mor{S\iii}{10}}\put(10,30){\ig{10}}\put(0,40){\ibraid}\put(0,50){\comult}\put(5,55){\comult}\put(10,60){\ig{5}}\put(20,50){\twist{-5}{5}}\put(26,30){=}\egr}
\put(35,0){\bgr{.4ex}{20}{85}\put(15,0){\ig{5}}\put(10,5){\mult}\put(0,4){\coumor{\va}{6}}\put(0,10){\isibrmor}\put(0,20){\sibrmor}\put(0,30){\ig{10}}\put(0,40){\isibrmor}\put(0,50){\ig{10}}\put(10,50){\mor{S}{10}}\put(0,60){\comult}\put(10,30){\sibrmor}\put(20,10){\ig{20}}\put(20,40){\ig{30}}\put(20,70){\twist{-5}{5}}\put(5,65){\mor{S\iii}{10}}\put(5,75){\comult}\put(10,80){\ig{5}}\put(26,40){=}\egr}
\put(70,0){\bgr{.4ex}{20}{85}\put(15,0){\ig{5}}\put(10,5){\mult}\put(0,4){\coumor{\va}{6}}\put(0,10){\isibrmor}\put(20,10){\ig{10}}\put(0,20){\ig{10}}\put(0,30){\sibrmor}\put(10,40){\sibrmor}\put(20,30){\ig{10}}\put(0,40){\ig{20}}\put(20,50){\ig{20}}\put(10,80){\ig{5}}\put(10,20){\isibrmor}\put(10,50){\mor{S}{10}}\put(0,60){\comult}\put(5,65){\mor{S\iii}{10}}\put(5,75){\comult}\put(20,70){\twist{-5}{5}}\put(26,40){=}\egr}
\put(105,0){\bgr{.4ex}{20}{85}\put(15,0){\ig{10}}\put(5,4){\coumor{\va}{6}}\put(5,10){\isibrmor}\put(5,20){\ig{5}}\put(15,20){\twist{5}{5}}\put(20,25){\ig{15}}\put(0,25){\mult}\put(0,30){\sibrmor}\put(0,40){\ig{20}}\put(0,60){\comult}\put(5,65){\mor{S\iii}{10}}\put(5,75){\comult}\put(20,70){\twist{-5}{5}}\put(20,50){\ig{20}}\put(10,40){\sibrmor}\put(10,50){\mor{S}{10}}\put(10,80){\ig{5}}\put(26,40){=}\egr}
\put(140,-5){\bgr{.4ex}{20}{85}\put(15,0){\ig{10}}\put(5,4){\coumor{\va}{6}}\put(5,10){\isibrmor}\put(15,20){\twist{5}{5}}\put(20,25){\mor{S}{10}}\put(0,20){\mult}\put(0,25){\sibrmor}\put(0,35){\ig{30}}\put(0,65){\comult}\put(5,70){\mor{S\iii}{10}}\put(5,80){\comult}\put(20,75){\twist{-5}{5}}\put(10,35){\braid}\put(10,45){\isibrmor}\put(10,55){\ibraid}\put(20,65){\ig{10}}\put(10,85){\ig{5}} \put(26,45){=}\egr}
\put(175,-5){\bgr{.4ex}{20}{85}\put(15,0){\ig{10}}\put(5,4){\coumor{\va}{6}}\put(5,10){\isibrmor}
\put(15,20){\twist{5}{5}}\put(20,25){\mor{S}{10}}\put(0,20){\mult}\put(0,25){\sibrmor}\put(0,35){\ig{10}}\put(10,35){\braid}
\put(0,45){\comult}\put(20,45){\twist{-5}{5}}\put(5,50){\isibrmor}\put(5,60){\ibraid}\put(5,70){\mor{S\iii}{10}}
\put(15,70){\ig{10}}\put(5,80){\comult}\put(10,85){\ig{5}}\put(26,45){=}\egr}
\cegr

\cbgr{75}{60}\eeee
\put(0,0){\bgr{.4ex}{20}{70}\put(20,0){\ig{10}}\put(10,4){\coumor{\va}{6}}\put(10,10){\isibrmor}\put(10,20){\twist{-5}{5}}\put(0,25){\mult}\put(0,30){\ig{10}}\put(10,30){\comult}\put(20,20){\mor{S}{10}}\put(15,35){\twist{-5}{5}}\put(0,40){\ibraid}\put(0,50){\mor{S}{10}}\put(10,50){\ig{10}}\put(0,60){\comult}\put(5,65){\ig{5}}\put(26,35){=}\egr}
\put(35,5){\bgr{.4ex}{20}{70}\put(15,0){\ig{10}}\put(5,4){\coumor{\va}{6}}\put(5,10){\isibrmor}\put(5,20){\ig{5}}\put(0,25){\mult}\put(10,30){\ig{10}}\put(0,30){\mor{S}{10}}\put(0,40){\twist{5}{5}}\put(5,45){\comult}\put(10,50){\mor{S\iii}{10}}\put(15,20){\twist{5}{10}}\put(20,30){\mor{S}{10}}\put(10,40){\comult}\put(26,30){=}\egr}
\put(70,12){\bgr{.4ex}{20}{70}\put(0,0){\umor{\eta}{20}}\put(0,25){\coumor{\va}{20}}\egr}\cegr

 There is another way of defining the opposite product: $m^{\op'}:=m\si^{-1}$. For this opposite product the antipode is give by
\bbas \overline{S}':=(S^{-1}\ot \va)\si\tau^{-1} \Delta.\eeas
We leave it to the reader to check this equation.

A central bialgebra $(B,\si)$  is said to be commutative if $m^\op=m$. If $B$ is a Hopf algebra then we have $\overline{S}=S$, since the antipode exists uniquely. Thus we have shown:
\begin{cor} Let $(B,\si)$ be a commutative central Hopf algebra, then the antipode of $B$ is invertible.\end{cor}
Note that the square of the antipode in not the identity morphism, unless the braiding is symmetric. The the next section we shall study the generalization of the notion commutative -- the coquasitriangular structures.

\section{Coquasitriangular structures on central bialgebras}\label{sec-cqt}
Since on a central bialgebra there exist two multiplication which both agree with the comultiplication, we can define coquasitriangular (CQT) structures which compare these two multiplications \cite{majid1}. The main results we obtain here is that, in a central Hopf algebra, which admits a CQT structure, the antipode is invertible and the square of the antipode can be given via the braiding, $\si$ and the CQT structure.

Let $(B,\si)$ be a central bialgebra. A coquasitriangular structure on $B$ is a morphism $r:B\ot B\lora I$ in $\C$ which obeys the following axioms:

\begin{description}

\item[CQT1] The two products can be compared by $r$:
\bb\label{eqcqt1}(m^{\op}\otimes r)\tau_{B,B}(\Delta\otimes \Delta)=(r\otimes m)\tau_{B,B}(\Delta\otimes\Delta),\ee

\item [CQT2] $r$ satisfies the following equations:
\bb\label{eqcqt21} r(m\ot B)=(r\ot r)(B\ot\si\ot B)(B\ot B\ot \Delta) ,\ee
\bb\label{eqcqt22} r(B\otimes m)=r(B\otimes r\otimes B)(\Delta\otimes B\ot B),\ee

\item[CQT3] There exists  a morphism $r^*$, subject to the following equations:
\bb\label{eqcqt3}\bbar{c}(r^*\si\otimes r)\Delta_{B\ot B} =(r\ot r^*\si)\Delta_{B\ot B} =\va\otimes\va .\eear\ee
\end{description}
\rref{eqcqt1}--\rref{eqcqt3} are depicted as follows:

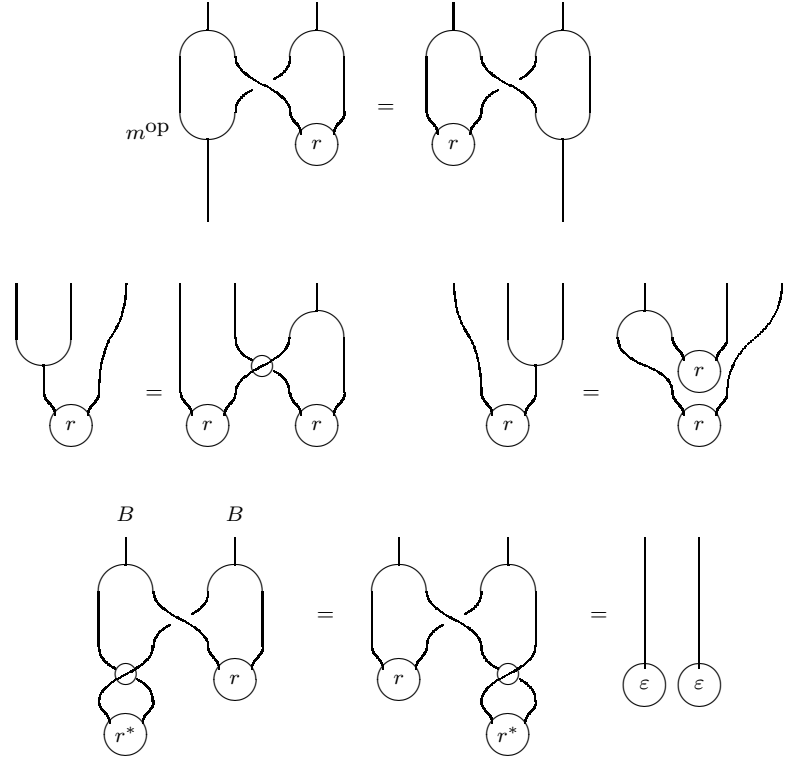
\cbgr{90}{40}
\put(50,0){\bgr{.4ex}{30}{40}\put(0,10){\rmor{r}}\put(0,20){\ig{10}}\put(0,30){\comult}\put(5,35){\ig{5}}\put(10,20){\braid}\put(25,0){\ig{15}}\put(20,15){\mult}\put(30,20){\ig{10}}\put(20,30){\comult}\put(25,35){\ig{5}}\egr}

\put(5,0){\bgr{.4ex}{30}{40}\put(-10,15){$m^{\op}$}\put(20,10){\rmor{r}}\put(0,20){\ig{10}}\put(0,30){\comult}\put(5,35){\ig{5}}\put(10,20){\braid}\put(5,0){\ig{15}}\put(0,15){\mult}\put(30,20){\ig{10}}\put(20,30){\comult}\put(25,35){\ig{5}}\put(36,20){=}\egr}\cegr

\cbgr{140}{30}
\put(0,0){\bgr{.4ex}{20}{35}\put(0,20){\mult}\put(0,25){\idgr{10}}\put(10,25){\idgr{10}}\put(5,15){\idgr{5}}\put(15,15){\twist{5}{20}}\put(5,5){\rmor{r}}\egr}
\put(25,15){=}

\put(30,5){\bgr{.4ex}{20}{35}
\put(20,0){\rmor{r}}
\put(10,10){\simor}
\put(10,20){\idgr{10}}
\put(0,0){\rmor{r}}
\put(20,20){\comult}
\put(25,25){\idgr{5}}
\put(0,10){\idgr{20}}
\put(30,10){\idgr{10}}
\egr}

\put(80,0){\bgr{.4ex}{20}{35}\put(10,20){\mult}\put(20,25){\idgr{10}}\put(10,25){\idgr{10}}\put(15,15){\idgr{5}}\put(5,15){\twist{-5}{20}}\put(5,5){\rmor{r}}
\egr}\put(105,15){=}

\put(110,0){\bgr{.4ex}{30}{35}\put(5,30){\idgr{5}}\put(20,25){\idgr{10}}\put(0,25){\comult}\put(0,25){\twist{10}{-10}}\put(10,5){\rmor{r}}\put(10,15){\rmor{r}}\put(20,15){\twist{10}{20}}\egr}
\cegr

\cbgr{110}{30}
\put(0,0){\bgr{.4ex}{30}{35}\put(0,15){\idgr{10}}\put(0,25){\comult}\put(5,30){\idgr{5}}\put(25,30){\idgr{5}}\put(30,15){\idgr{10}}\put(0,5){\simor}\put(0,-5){\rmor{r^*}}\put(20,5){\rmor{r}}\put(10,15){\braid}\put(20,25){\comult}\put(05,35){\objo{B}}\put(25,35){\objo{B}}\put(40,20){=}\egr}
\put(50,0){\bgr{.4ex}{30}{35}\put(0,15){\idgr{10}}\put(0,25){\comult}\put(5,30){\idgr{5}}\put(25,30){\idgr{5}}\put(30,15){\idgr{10}}\put(20,5){\simor}\put(20,-5){\rmor{r^*}}\put(0,5){\rmor{r}}\put(10,15){\braid}\put(20,25){\comult}
\put(40,20){=}\egr}
\put(100,-5){\bgr{.4ex}{30}{35}\put(0,10){\coumor{\va}{30}}\put(10,10){\coumor{\va}{30}}
\egr}
\cegr

\noindent{\bf Remarks.} Axiom CQT3 implies that $r$ is $*$-invertible, $r^{-1}=r^*\si^{-1}_{B,B}$. We use $r*$ instead of $r^{-1}$ just because the computation with $r^*$ is simpler.  Our definition of CQT structure slightly differs from Majid's definition, in which \rref{eqcqt21} reads:
$r (m^{\op}\otimes B)=r(B\otimes R\otimes B)(B\otimes B\otimes \Delta)$. Note however, that in the example, given by Majid, of reconstructed bialgebras \cite{majid2}, the two definitions are equivalent. In fact, the two definitions will be equivalent if $r$ and $\si$ commute, i.e., if 
\bb\label{eqrsi}(B\ot r)\si_{B,B\ot B}=r\ot B.\ee 
In our setting, this equation holds, for example, if $r$ is a comodule morphism. However, $r$ is in general not a comodule morphism, except when $r=\va\ot \va$. For a reconstructed bialgebra $(B,\si)$, $(B,\si)$ is also object of $\Zz(\C)$, hence, \rref{eqrsi} holds.

The following equation follows immediately from \rref{eqcqt1} and \rref{eqcqt3}:
\bb\label{rsm}(r^*\si\ot m\si)\Delta_{B\ot B}=(m\ot r^*\si)\Delta_{B\ot B}.\ee
In graphical notations \rref{rsm} has the form:

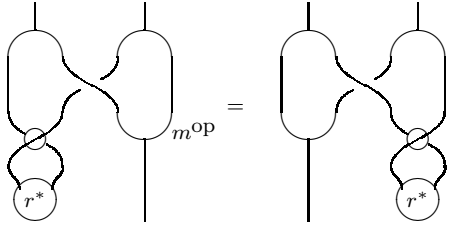
\cbgr{90}{40}
\put(0,0){\bgr{.4ex}{30}{40}\put(0,10){\simor}\put(0,0){\rmor{r^*}}\put(0,20){\ig{10}}\put(0,30){\comult}\put(5,35){\ig{5}}\put(10,20){\braid}\put(25,0){\ig{15}}\put(20,15){\mult}\put(30,20){\ig{10}}\put(20,30){\comult}\put(25,35){\ig{5}}\put(40,20){=}\put(30,15){$m^{\op}$}\egr}

\put(50,0){\bgr{.4ex}{30}{40}\put(20,10){\simor}\put(20,0){\rmor{r^*}}\put(0,20){\ig{10}}\put(0,30){\comult}\put(5,35){\ig{5}}\put(10,20){\braid}\put(5,0){\ig{15}}\put(0,15){\mult}\put(30,20){\ig{10}}\put(20,30){\comult}\put(25,35){\ig{5}}\egr}\cegr

\begin{lem}\label{lemrrs} Let  $B$ be a central, CQT Hopf algebra then $r,r^*$ obey the following equations:
\bb\label{eqSrs1} r(S\otimes\id)=r^*\tau^{-1} ,\ee \bb \label{eqSrs2}
r^*(S\otimes\id)=r \si^{-1},\ee

\cbgr{80}{15}\put(0,-5){\bgr{.4ex}{10}{10}
\put(0,0){\rmor{r}}\put(0,10){\mor{S}{10}}\put(20,10){\ibraid}\put(10,10){\idgr{10}}\put(15,15){=}\put(20,0){\rmor{r^*}}
\put(50,0){\rmor{r^*}}\put(50,10){\mor{S}{10}}\put(60,10){\idgr{10}}\put(65,10){=}
\put(70,0){\rmor{r}}\put(70,10){\isibrmor}
\put(0,20){\objo{B}}\put(10,20){\objo{B}}\put(20,20){\objo{B}}\put(30,20){\objo{B}}\put(50,20){\objo{B}}\put(60,20){\objo{B}}\put(70,20){\objo{B}}\put(80,20){\objo{B}}\egr}
\cegr

\bb\label{eqmrs1} r^*(m\otimes B)=(r^*\ot r^*)(B\otimes \tau\otimes B)(B\otimes  B\otimes \Delta),\ee \bb \label{eqmrs2}
r^*(B\otimes m )=(r^*\ot r^*)(B\otimes \si^{-1}\otimes B)(\tau^{-1}\Delta\otimes B\ot B).
\ee

\cbgr{150}{30}
\put(0,-5){\bgr{.4ex}{20}{35}\put(0,20){\mult}\put(0,25){\idgr{10}}\put(10,25){\idgr{10}}\put(5,15){\idgr{5}}\put(15,15){\twist{5}{20}}\put(5,5){\rmor{r^*}}\put(0,35){\objo{B}}\put(10,35){\objo{B}}\put(20,35){\objo{B}}\egr}\put(25,15){=}
\put(33,-5){\bgr{.4ex}{30}{35}\put(0,15){\idgr{20}}\put(10,25){\idgr{10}}\put(25,30){\idgr{5}}\put(30,15){\idgr{10}}\put(0,5){\rmor{r^*}}\put(20,5){\rmor{r^*}}\put(10,15){\braid}\put(20,25){\comult}\put(0,35){\objo{B}}\put(10,35){\objo{B}}\put(25,35){\objo{B}}\egr}
\put(80,-5){\bgr{.4ex}{20}{35}\put(10,20){\mult}\put(20,25){\idgr{10}}\put(10,25){\idgr{10}}\put(15,15){\idgr{5}}\put(5,15){\twist{-5}{10}}\put(5,5){\rmor{r^*}}\put(0,25){\ig{10}}\put(0,35){\objo{B}}\put(10,35){\objo{B}}\put(20,35){\objo{B}}\egr}\put(105,15){=}
\put(113,0){\bgr{.3ex}{20}{35}\put(0,0){\rmor{r^*}}\put(0,10){\ig{10}}\put(0,30){\comult}\put(5,35){\idgr{5}}\put(0,20){\ibraid}\put(10,10){\isibrmor}\put(20,0){\rmor{r^*}}\put(30,10){\idgr{30}}\put(20,20){\idgr{20}}\put(05,40){\objo{B}}\put(20,40){\objo{B}}\put(30,40){\objo{B}}\egr}\cegr
\end{lem}
 
\proof Using \rref{eqcqt21} we have
\cbgr{100}{55}
\put(0,10){\bgr{.4ex}{30}{50}\put(0,-10){\rmor{r}}\put(0,0){\mor{S}{10}}\put(0,10){\braid}\put(0,20){\sibrmor}\put(0,30){\ig{10}}\put(0,40){\comult}\put(5,45){\ig{5}}\put(10,0){\ig{10}}\put(10,30){\braid}\put(20,20){\rmor{r}}\put(30,30){\ig{10}}\put(20,40){\comult}\put(25,45){\ig{5}}\put(36,20){=}\egr}
\put(45,15){\bgr{.4ex}{30}{50}\put(0,0){\rmor{r}}\put(0,10){\mor{S}{10}}\put(0,20){\comult}\put(10,10){\sibrmor}\put(20,20){\comult}\put(20,0){\rmor{r}}\put(5,25){\ig{5}}\put(25,25){\ig{5}}\put(30,10){\ig{10}}\put(36,15){=}\egr}
\put(90,5){\bgr{.4ex}{30}{50}\put(0,10){\coumor{\va}{30}}\put(10,10){\coumor{\va}{30}}\egr}
\cegr
Thus $r(S\ot B)\tau\si$ is the left $*$-inverse of $r$ in $\Mor(B\ot B, I)$. Since $r$ is $*$-invertible, \rref{eqSrs1} follows.
For \rref{eqSrs2} we have:

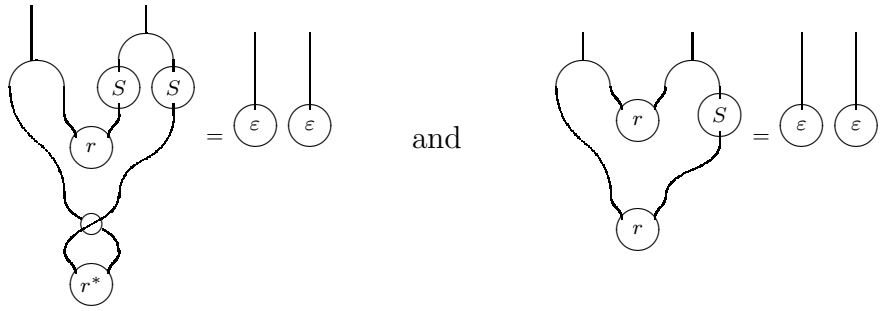
\cbgr{155}{55}
\put(0,0){\bgr{.4ex}{30}{50}\put(10,0){\rmor{r^*}}\put(10,10){\sibrmor}\put(10,20){\twist{-10}{20}}\put(0,40){\comult}\put(4,45){\ig{10}}\put(10,35){\ig{5}}\put(10,25){\rmor{r}}\put(20,35){\mor{S}{10}}\put(20,45){\comult}\put(25,50){\ig{5}}\put(30,35){\mor{S}{10}}\put(20,20){\twist{10}{15}}\put(36,30){=}\egr}
\put(45,0){\bgr{.4ex}{30}{50}\put(0,30){\coumor{\va}{20}}\put(10,30){\coumor{\va}{20}}\egr}\put(75,30){\nml and}
\put(100,10){\bgr{.4ex}{30}{50}\put(10,0){\rmor{r}}\put(10,10){\twist{-10}{20}}\put(0,30){\comult}\put(5,35){\ig{5}}\put(20,10){\twist{10}{10}}\put(30,20){\mor{S}{10}}\put(10,20){\rmor{r}}\put(20,30){\comult}\put(25,35){\ig{5}}\put(36,20){=}\egr}
\put(145,0){\bgr{.4ex}{30}{50}\put(0,30){\coumor{\va}{20}}\put(10,30){\coumor{\va}{20}}\egr}\cegr
The first equation follows from \rref{eqcqt3}, the second one follows from \rref{eqcqt22}.
 To prove \rref{eqmrs1} we use  \rref{eqSrs1}, \rref{eqSrs2} and \rref{eqcqt22}:

\cbgr{128}{50}
\put(0,10){\bgr{.4ex}{30}{50}\put(0,0){\rmor{r^*}}\put(0,10){\ig{20}}\put(10,10){\braid}\put(10,20){\ig{10}}\put(20,0){\rmor{r^*}}\put(20,20){\comult}\put(30,10){\ig{10}}\put(25,25){\ig{5}}\put(33,10){=}\egr}
\put(40,10){\bgr{.4ex}{30}{50}\put(0,-10){\rmor{r}}\put(0,0){\braid}\put(0,10){\ig{30}}\put(10,10){\mor{S}{10}}\put(10,20){\braid}\put(20,0){\rmor{r}}\put(20,10){\braid}\put(30,20){\mor{S}{10}}\put(10,30){\ig{10}}\put(20,30){\comult}\put(25,35){\ig{5}}\put(36,10){=}\egr}
\put(80,5){\bgr{.4ex}{30}{50}\put(0,30){\ig{5}}\put(5,0){\rmor{r}}\put(5,10){\braid}\put(5,20){\ig{5}}\put(0,25){\mult}\put(10,30){\ig{5}}\put(15,20){\mor{S}{10}}\put(15,30){\ig{5}}\put(21,15){=}\egr}
\put(105,15){\bgr{.4ex}{30}{50}\put(0,20){\ig{5}}\put(5,10){\ig{5}}\put(10,20){\ig{5}}\put(0,15){\mult}\put(5,0){\rmor{r^*}}\put(15,10){\ig{15}}\egr}\cegr

To prove \rref{eqmrs2} we use  \rref{sigmaantipode}  and \rref{eqcqt21}:

\cbgr{120}{70}
\put(0,20){\bgr{.4ex}{40}{40}\put(0,0){\rmor{r^*}}\put(20,0){\rmor{r^*}}\put(0,10){\ig{10}}\put(0,30){\comult}\put(5,35){\ig{5}}\put(10,10){\isibrmor}\put(0,20){\ibraid}\put(20,20){\ig{20}} \put(30,10){\ig{30}}\egr}\put(36,40){=}
\put(45,30){\bgr{.4ex}{30}{50}\put(0,-20){\rmor{r}}\put(20,-20){\rmor{r}}\put(0,10){\ig{10}}\put(0,30){\comult}\put(5,35){\ig{5}}\put(10,10){\isibrmor}\put(0,20){\ibraid}\put(20,20){\ig{20}} \put(30,10){\ig{30}}\put(0,0){\braid}\put(20,0){\braid}\put(0,-10){\mor{S}{10}}\put(20,-10){\mor{S}{10}}\put(10,-10){\ig{10}}\put(30,-10){\ig{10}}\egr}\put(81,40){=}

\put(90,0){\bgr{.4ex}{30}{50}\put(0,0){\rmor{r}}\put(20,0){\rmor{r}}\put(0,10){\braid}\put(20,10){\braid}\put(0,20){\ig{30}}\put(10,20){\ibraid}\put(30,20){\mor{S}{50}}\put(10,30){\sibrmor}\put(20,50){\mor{S}{20}}\put(10,40){\braid}\put(0,50){\ibraid}\put(0,60){\comult}\put(5,65){\ig{5}}\egr}\put(126,40){=}
\cegr

\cbgr{90}{40}
\put(0,-5){\bgr{.4ex}{30}{50}\put(0,0){\rmor{r}}\put(20,0){\rmor{r}}\put(0,10){\ig{10}}\put(10,10){\sibrmor}\put(0,20){\braid}\put(20,20){\comult}\put(30,10){\ig{10}}\put(25,25){\twist{-5}{5}}\put(10,30){\braid}\put(0,30){\ig{10}}\put(0,40){\braid}\put(0,50){\ig{10}}\put(10,50){\mor{S}{10}}\put(20,40){\mor{S}{20}}\egr}\put(36,25){=}

\put(45,0){\bgr{.4ex}{20}{50}\put(5,-5){\rmor{r}}\put(0,5){\mult}\put(0,10){\braid}\put(0,20){\ig{10}}\put(0,30){\braid}\put(0,40){\ig{10}}\put(15,5){\twist{5}{15}}\put(10,20){\braid}\put(20,30){\mor{S}{20}}\put(10,40){\mor{S}{10}}\egr}\put(76,25){=}
\put(80,10){\bgr{.4ex}{10}{50}\put(0,0){\rmor{r^*}}\put(0,10){\ig{20}}\put(5,20){\ig{10}}\put(15,20){\ig{10}}\put(5,15){\mult}\put(10,10){\ig{5}}\egr}
\eeee\cegr

Composing both sides of \rref{eqcqt21} and \rref{eqcqt22}  with $(B\ot\eta \va\ot B)(B\ot \Delta)$ from the right and using \rref{eqz4}, we have:
\bbas r=r(B\ot r(\eta\va\ot B)\ot B)\Delta_{B\ot B},\\
r=r(B\ot r(B\ot\eta\va)\ot B)\Delta_{B\ot B}.\eeas
That is, $r=r*r(B\ot \eta\va)=r(\eta\va\ot B)*r$ in $\Mor(B\ot B,I)$. According to CQT3, $r$ is $*$-invertible,  hence 
\bb\label{eqreta}r(\eta\ot B)=r(B\ot \eta)=\va\ot \va.\ee

   \rref{eqSrs1}, \rref{eqSrs2} and \rref{eqreta} imply:
\bb\label{eqrseta}r^*(\eta\ot B)=r^*(B\ot \eta)=\va\ot\va .\ee

Let $(B,\si,r)$ be a CQT central Hopf algebra. We define  morphisms $u,u^*:B\lora I$ as as follows:
\bb\label{eqdefu}\bbar{cc} u= r \si^{-1}(S\otimes B),&
 u^*=r^*(B\otimes S).\eear\ee
According to   \rref{eqSrs1}, \rref{eqSrs2}, we have:
\bb\label{eqdefu1}\bbar{cc} u=r^*(S^2\ot B),& u^*=r\si^{-1}_{B,B}(S^*\ot B).\eear\ee
In graphical  notations the equations above are depicted as:

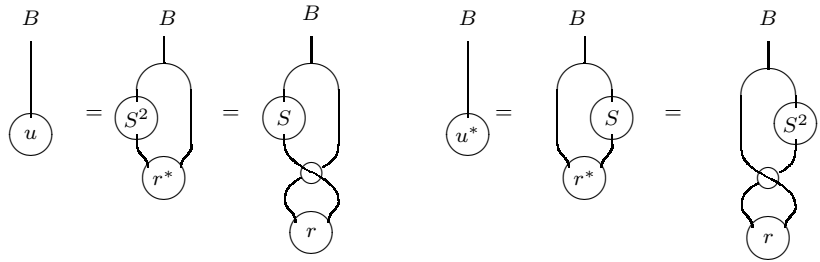
\cbgr{130}{40}
\put(0,-5){\bgr{.4ex}{20}{20}\put(0,20){\coumor{u}{20}}\put(10,26){=}\put(13,-5){\bgr{.4ex}{15}{30}\put(10,40){\idgr{5}}\put(5,35){\comult}\put(15,25){\idgr{10}}\put(5,25){\mor{S^2}{10}}\put(5,15){\rmor{r^*}}\egr}
\put(0,40){\objo{B}}\put(25,40){\objo{B}}\put(35,26){=}
\put(40,-10){\bgr{.4ex}{20}{20}\put(5,40){\comult}\put(15,30){\idgr{10}}\put(5,30){\mor{S}{10}}\put(5,20){\isibrmor}\put(5,10){\rmor{r}}\put(10,45){\idgr{5}}\put(10,50){\objo{B}}\egr}
\put(80,20){\coumor{u^*}{20}}\put(85,26){=}
\put(90,-5){\bgr{.4ex}{15}{30}\put(10,40){\idgr{5}}\put(5,35){\comult}\put(5,25){\idgr{10}}\put(15,25){\mor{S}{10}}\put(5,15){\rmor{r^*}}\egr}
\put(80,40){\objo{B}}\put(100,40){\objo{B}}\egr}
\put(125,-10){\bgr{.4ex}{20}{20}\put(10,40){\idgr{5}}\put(5,35){\comult}\put(5,25){\idgr{10}}\put(15,25){\mor{S^2}{10}}\put(5,5){\rmor{r}}\put(5,15){\isimor}\put(10,45){\objo{B}}\put(-9,31){=}\egr}

\cegr

\begin{lem} \label{lemuu} $u$ and $u^*$ satisfy the following  equations:
\bb\label{equu}
 (u\otimes B)\Delta=(u\otimes S^2)\si^{-1}\Delta,\ee
\bb \label{equsus}(u^*\ot S^2)\Delta= (B\ot u^*)\si\tau^{-1}\Delta,\ee
or in graphical notations:

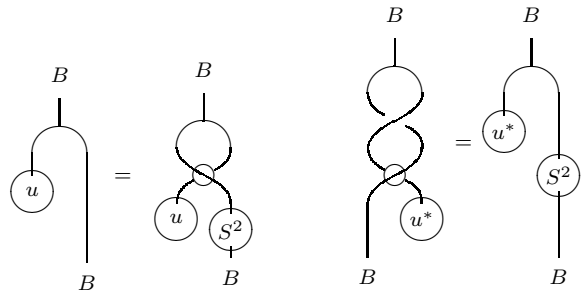
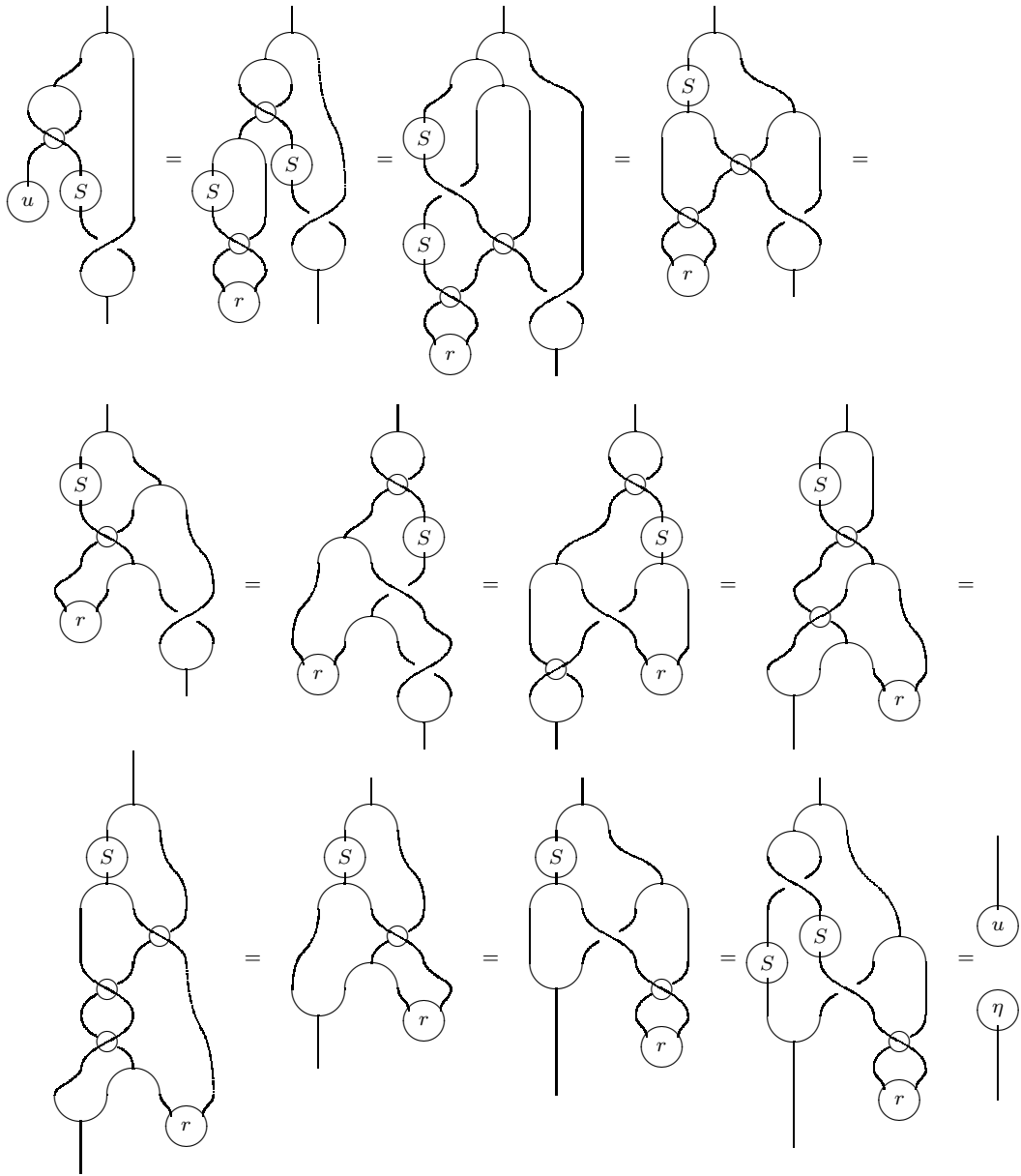
\cbgr{100}{40}
\put(0,10){\coumor{u}{10}}\put(0,20){\comult}\put(5,25){\idgr{5}}\put(10,0){\idgr{20}}
\put(10,0){\obju{B}}\put(5,30){\objo{B}}
\put(15,15){=}
\put(20,-5){\bgr{.4ex}{10}{40}
\put(5,25){\comult}\put(10,30){\idgr{5}}
\put(5,15){\isibrmor}\put(5,9){\coumor{u}{6}}\put(15,5){\mor{S^2}{10}}
\put(15,05){\obju{B}}\put(10,35){\objo{B}}\egr}

\put(60,0){\bgr{.4ex}{10}{50}\put(0,0){\ig{10}}\put(0,10){\sibrmor}\put(0,20){\ibraid}\put(0,30){\comult}\put(5,35){\ig{5}}\put(10,4){\coumor{u^*}{6}}\put(16,20){=}\put(0,0){\obju{B}}\put(5,40){\objo{B}}\egr}
\put(85,0){\bgr{.4ex}{10}{50}\put(10,0){\mor{S^2}{30}}\put(0,30){\comult}\put(5,35){\ig{5}}\put(0,20){\coumor{u^*}{10}}\put(10,0){\obju{B}}\put(5,40){\objo{B}}\egr}\cegr
\end{lem}
\proof 
We have:

\cbgr{190}{70}
\put(0,10){\bgr{.4ex}{20}{60}\put(15,0){\ig{5}}\put(10,5){\mult}\put(10,10){\ibraid}\put(10,20){\mor{S}{10}}\put(0,20){\coumor{u}{10}}\put(0,30){\isibrmor}\put(0,40){\comult}\put(5,45){\twist{5}{5}}\put(10,50){\comult}\put(15,55){\ig{5}}\put(20,20){\ig{30}}\put(26,30){=}\egr}
\put(35,10){\bgr{.4ex}{20}{60}\put(20,0){\ig{10}}\put(0,0){\rmor{r}}\put(10,20){\ig{10}}\put(0,10){\isibrmor}\put(0,20){\mor{S}{10}}\put(0,30){\comult}\put(5,35){\isibrmor}\put(5,45){\comult}\put(10,50){\comult}\put(15,55){\ig{5}}\put(15,10){\mult}\put(15,15){\ibraid}\put(15,25){\mor{S}{10}}\put(25,25){\twist{-5}{25}}\put(31,30){=}\egr}
\put(75,0){\bgr{.4ex}{20}{60}\put(0,0){\rmor{r}}\put(0,10){\isibrmor}\put(0,20){\mor{S}{10}}\put(0,30){\braid}\put(0,40){\mor{S}{10}}\put(10,40){\ig{10}}\put(0,50){\twist{5}{5}}\put(5,55){\comult}\put(10,60){\comult}\put(15,65){\ig{5}}\put(10,50){\comult}\put(10,20){\isibrmor}\put(20,30){\ig{20}}\put(25,0){\ig{5}}\put(20,5){\mult}\put(20,10){\ibraid}\put(30,20){\ig{30}}\put(30,50){\twist{-10}{10}}\put(36,40){=}\egr}
\put(120,15){\bgr{.4ex}{20}{60}\put(0,0){\rmor{r}}\put(0,10){\isibrmor}\put(0,20){\ig{10}}\put(0,30){\comult}\put(10,20){\isibrmor}\put(5,35){\mor{S}{10}}\put(5,45){\comult}\put(10,50){\ig{5}}\put(25,0){\ig{5}}\put(20,5){\mult}\put(20,10){\braid}\put(30,20){\ig{10}}\put(20,30){\comult}\put(25,35){\twist{-10}{10}}\put(36,25){=}\egr}\cegr
\cbgr{180}{55}
\put(0,10){\bgr{.4ex}{20}{60}\put(0,10){\rmor{r}}\put(0,20){\twist{5}{5}}\put(5,25){\isibrmor}\put(5,35){\mor{S}{10}}\put(5,45){\comult}\put(10,50){\ig{5}}\put(25,0){\ig{5}}\put(20,5){\mult}\put(20,10){\ibraid}\put(10,20){\comult}\put(30,20){\twist{-5}{15}}\put(15,35){\comult}\put(20,40){\twist{-5}{5}}\put(36,20){=}\egr}
\put(45,0){\bgr{.4ex}{20}{60}\put(25,0){\ig{5}}\put(20,5){\mult}\put(20,10){\ibraid}\put(10,20){\comult}\put(0,10){\rmor{r}}\put(0,20){\twist{5}{15}}\put(5,35){\comult}\put(10,40){\twist{5}{5}}\put(15,25){\braid}\put(30,20){\twist{-5}{5}}\put(25,35){\mor{S}{10}}\put(15,45){\isibrmor}\put(15,55){\comult}\put(20,60){\ig{5}}\put(36,30){=}\egr}
\put(90,0){\bgr{.4ex}{20}{60}\put(5,0){\ig{5}}\put(0,5){\mult}\put(0,10){\sibrmor}\put(0,20){\ig{10}}\put(0,30){\comult}\put(10,20){\braid}\put(20,10){\rmor{r}}\put(30,20){\ig{10}}\put(20,30){\comult}\put(25,35){\mor{S}{10}}\put(15,45){\isibrmor}\put(15,55){\comult}\put(20,60){\ig{5}}\put(5,35){\twist{10}{10}}\put(36,30){=}\egr}
\put(135,0){\bgr{.4ex}{20}{60}\put(5,0){\ig{10}}\put(0,10){\mult}\put(0,15){\twist{5}{5}}\put(10,15){\comult}\put(15,30){\comult}\put(20,5){\rmor{r}}\put(30,15){\twist{-5}{15}}\put(5,20){\sibrmor}\put(5,30){\twist{5}{5}}\put(10,35){\isibrmor}\put(10,45){\mor{S}{10}}\put(20,45){\ig{10}}\put(10,55){\comult}\put(15,60){\ig{5}}\put(36,30){=}\egr}\cegr
\cbgr{180}{65}
\put(0,0){\bgr{.4ex}{30}{60}\put(5,0){\ig{10}}\put(0,10){\mult}\put(0,15){\twist{5}{5}}\put(10,15){\comult}\put(20,5){\rmor{r}}\put(5,20){\sibrmor}\put(5,30){\isibrmor}\put(30,15){\twist{-5}{25}}\put(5,40){\ig{10}}\put(5,50){\comult}\put(10,55){\mor{S}{10}}\put(15,40){\isibrmor}\put(25,50){\twist{-5}{15}}\put(10,65){\comult}\put(15,70){\ig{10}}\put(36,40){=}\egr}
\put(45,20){\bgr{.4ex}{30}{60}\put(5,0){\ig{10}}\put(0,10){\mult}\put(0,15){\twist{5}{15}}\put(10,15){\comult}\put(20,5){\rmor{r}}\put(15,20){\isibrmor}\put(5,30){\comult}\put(10,35){\mor{S}{10}}\put(10,45){\comult}\put(15,50){\ig{5}}\put(30,15){\twist{-5}{5}}\put(25,30){\twist{-5}{15}}\put(36,20){=}\egr}
\put(90,15){\bgr{.4ex}{30}{60}\put(5,0){\ig{20}}\put(0,20){\mult}\put(0,25){\ig{10}}\put(0,35){\comult}\put(5,40){\mor{S}{10}}\put(5,50){\comult}\put(10,55){\ig{5}}\put(20,5){\rmor{r}}\put(20,15){\isibrmor}\put(30,25){\ig{10}}\put(10,25){\braid}\put(20,35){\comult}\put(25,40){\twist{-10}{10}}\put(36,25){=}\egr}
\put(135,5){\bgr{.4ex}{30}{60}\put(5,0){\ig{20}}\put(0,20){\mult}\put(0,25){\mor{S}{20}}\put(0,45){\braid}\put(10,25){\braid}\put(10,35){\mor{S}{10}}\put(0,55){\comult}\put(20,5){\rmor{r}}\put(20,15){\isibrmor}\put(30,25){\ig{10}}\put(20,35){\comult}\put(25,40){\twist{-10}{20}}\put(5,60){\comult}\put(10,65){\ig{5}}\put(36,35){=}\egr}
\put(180,45){\coumor{u}{20}} \put(180,15){\umor{\eta}{20}}  
\cegr

In the second equation of  the second row we use axiom CQT1. 
Thus composing both sides of the equation with $S$ we get:

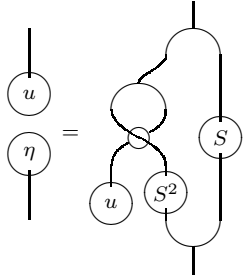
\cbgr{50}{50} \put(0,10){\bgr{.4ex}{20}{60}\put(0,20){\coumor{u}{15}}\put(0,0){\umor{\eta}{15}}\put(6,15){=}\egr}
\put(15,-10){\bgr{.4ex}{20}{60}\put(15,10){\ig{5}}\put(10,15){\mult}\put(10,20){\mor{S^2}{10}}\put(0,20){\coumor{u}{10}}\put(0,30){\isibrmor}\put(0,40){\comult}\put(5,45){\twist{5}{5}}\put(10,50){\comult}\put(15,55){\ig{5}}\put(20,20){\mor{S}{30}}\egr}\cegr
which proves \rref{equu}. For proving \rref{equsus} we have:

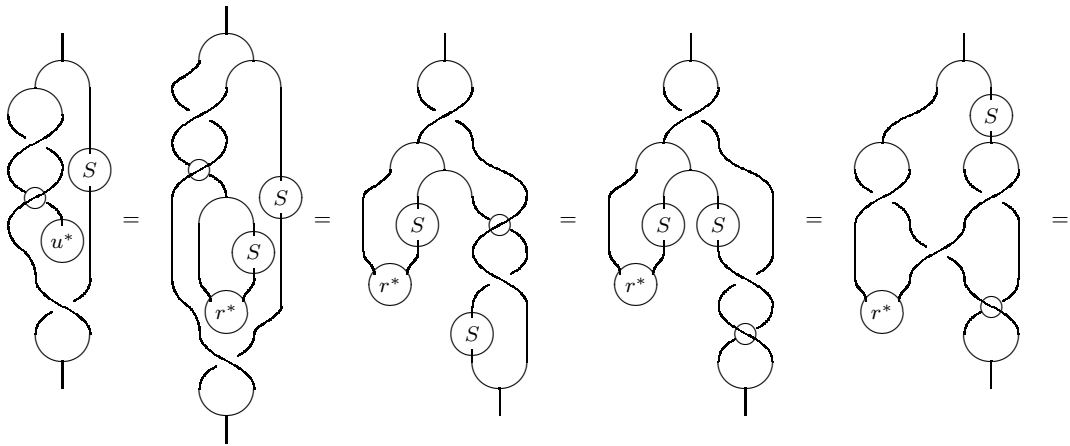
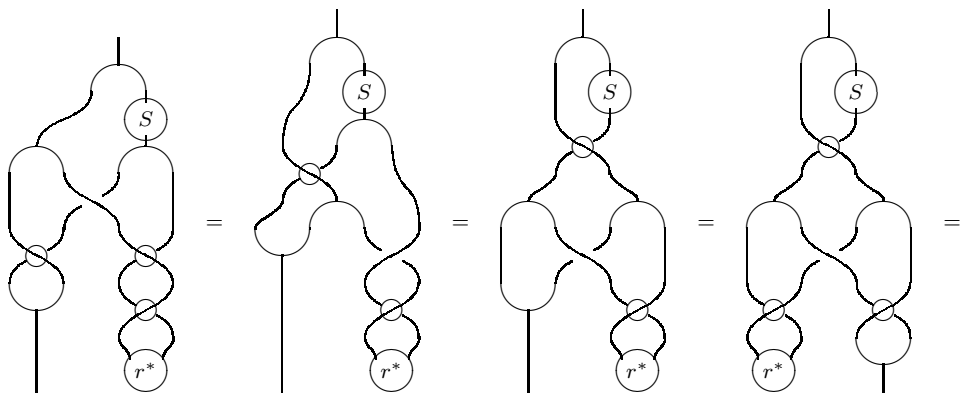
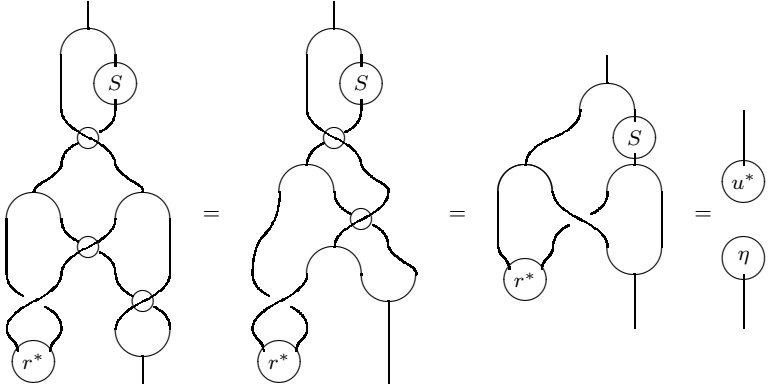
\cbgr{190}{75}\put(0,10){\bgr{.4ex}{30}{70}\put(10,60){\ig{5}}\put(5,55){\comult}\put(0,50){\comult}\put(0,40){\ibraid}\put(0,30){\sibrmor}\put(5,20){\twist{-5}{10}}\put(5,10){\braid}\put(5,5){\mult}\put(10,0){\ig{5}}\put(10,24){\coumor{u^*}{6}}\put(15,20){\ig{5}}\put(15,25){\mor{S}{30}}\put(21,30){=}\egr}
\put(30,0){\bgr{.4ex}{30}{70}\put(10,0){\ig{5}}\put(5,5){\mult}\put(5,10){\braid}\put(5,20){\twist{-5}{10}}\put(0,30){\ig{15}}\put(0,45){\sibrmor}\put(0,55){\ibraid}\put(0,65){\twist{5}{5}}\put(5,70){\comult}\put(10,75){\ig{5}}\put(10,65){\comult}\put(15,20){\twist{5}{5}}\put(20,25){\mor{S}{40}}\put(5,20){\rmor{r^*}}\put(5,30){\ig{10}}\put(15,30){\mor{S}{10}}\put(5,40){\comult}\put(26,40){=}\egr}
\put(65,5){\bgr{.4ex}{30}{70}\put(0,40){\twist{5}{5}}\put(5,45){\comult}\put(10,50){\ibraid}\put(10,60){\comult}\put(15,65){\ig{5}}\put(10,40){\comult}\put(30,40){\twist{-10}{10}}  \put(0,20){\rmor{r^*}}\put(0,30){\ig{10}}\put(10,30){\mor{S}{10}}\put(25,0){\ig{5}}\put(20,5){\mult}\put(20,10){\mor{S}{10}}\put(30,10){\ig{10}}\put(20,20){\braid}\put(20,30){\sibrmor}\put(36,35){=}\egr}
\put(110,5){\bgr{.4ex}{30}{70}\put(0,40){\twist{5}{5}}\put(5,45){\comult}\put(10,50){\ibraid}\put(10,60){\comult}\put(15,65){\ig{5}}\put(10,40){\comult}\put(30,40){\twist{-10}{10}}\put(0,20){\rmor{r^*}}\put(0,30){\ig{10}}\put(10,30){\mor{S}{10}}\put(20,30){\mor{S}{10}}\put(25,0){\ig{5}}\put(20,5){\mult}\put(20,10){\isibrmor}\put(20,20){\braid}\put(30,30){\ig{10}}\put(36,35){=}\egr}
\put(155,10){\bgr{.4ex}{30}{70}\put(0,40){\comult}\put(20,40){\comult}\put(5,45){\twist{10}{10}}\put(25,45){\mor{S}{10}}\put(15,55){\comult}\put(20,60){\ig{5}} \put(0,10){\rmor{r^*}}\put(0,20){\ig{10}}\put(0,30){\ibraid}\put(10,20){\ibraid}\put(20,30){\ibraid}\put(20,10){\isibrmor}\put(30,20){\ig{10}}\put(25,0){\ig{5}}\put(20,5){\mult}\put(36,30){=}\egr}\cegr
\cbgr{180}{65}
\put(0,0){\bgr{.4ex}{30}{70}\put(0,40){\comult}\put(20,40){\comult}\put(5,45){\twist{10}{10}}\put(25,45){\mor{S}{10}}\put(15,55){\comult}\put(20,60){\ig{5}} \put(5,0){\ig{15}}\put(0,15){\mult}\put(0,20){\isibrmor}\put(0,30){\ig{10}}\put(20,0){\rmor{r^*}}\put(20,10){\sibrmor}\put(20,20){\isibrmor}\put(10,30){\braid}\put(30,30){\ig{10}}\put(36,30){=}\egr}
\put(45,0){\bgr{.4ex}{30}{70}\put(5,0){\ig{25}}\put(0,25){\mult}\put(0,30){\twist{5}{5}}\put(10,30){\comult}\put(5,35){\isibrmor}\put(5,45){\twist{5}{15}}\put(10,60){\comult}\put(15,65){\ig{5}}\put(20,50){\mor{S}{10}}\put(15,45){\comult}\put(20,0){\rmor{r^*}}\put(20,10){\sibrmor}\put(20,20){\ibraid}\put(30,30){\twist{-5}{15}}\put(36,30){=}\egr}
\put(90,0){\bgr{.4ex}{30}{70}\put(0,20){\ig{10}}\put(0,30){\comult}\put(10,20){\braid}\put(30,20){\ig{10}}\put(20,30){\comult}\put(5,35){\twist{5}{5}}\put(25,35){\twist{-5}{5}}\put(10,40){\isibrmor}\put(10,50){\ig{10}}\put(10,60){\comult}\put(15,65){\ig{5}}\put(20,50){\mor{S}{10}}         \put(5,0){\ig{15}}\put(0,15){\mult}\put(20,0){\rmor{r^*}}\put(20,10){\sibrmor}\put(36,30){=}\egr}
\put(135,0){\bgr{.4ex}{30}{70}\put(0,20){\ig{10}}\put(0,30){\comult}\put(10,20){\braid}\put(30,20){\ig{10}}\put(20,30){\comult}\put(5,35){\twist{5}{5}}\put(25,35){\twist{-5}{5}}\put(10,40){\isibrmor}\put(10,50){\ig{10}}\put(10,60){\comult}\put(15,65){\ig{5}}\put(20,50){\mor{S}{10}}      \put(0,0){\rmor{r^*}}\put(0,10){\sibrmor}\put(25,0){\ig{5}}\put(20,5){\mult}\put(20,10){\sibrmor}\put(36,30){=}\egr}
\cegr 
\cbgr{180}{60}
\put(0,0){\bgr{.4ex}{30}{70}\put(0,20){\ig{10}}\put(0,30){\comult}\put(10,20){\sibrmor}\put(30,20){\ig{10}}\put(20,30){\comult}\put(5,35){\twist{5}{5}}\put(25,35){\twist{-5}{5}}\put(10,40){\isibrmor}\put(10,50){\ig{10}}\put(10,60){\comult}\put(15,65){\ig{5}}\put(20,50){\mor{S}{10}}      \put(0,0){\rmor{r^*}}\put(0,10){\ibraid}\put(25,0){\ig{5}}\put(20,5){\mult}\put(20,10){\sibrmor}\put(36,30){=}\egr}
\put(45,0){\bgr{.4ex}{30}{70} \put(0,0){\rmor{r^*}}\put(0,10){\ibraid}\put(25,0){\ig{15}}\put(20,15){\mult}\put(0,20){\twist{5}{15}}\put(5,35){\comult}\put(10,40){\isibrmor}\put(10,50){\ig{10}}\put(20,50){\mor{S}{10}}\put(10,60){\comult}\put(15,65){\ig{5}}\put(10,20){\comult}\put(15,25){\sibrmor}\put(30,20){\twist{-5}{5}}\put(25,35){\twist{-5}{5}}\put(36,30){=}\egr}
\put(90,10){\bgr{.4ex}{30}{70} \put(0,05){\rmor{r^*}}\put(0,15){\ig{10}}\put(0,25){\comult}\put(10,15){\braid}\put(25,0){\ig{10}}\put(20,10){\mult}\put(30,15){\ig{10}}\put(0,25){\comult}\put(20,25){\comult}\put(5,30){\twist{10}{10}}\put(15,40){\comult}\put(20,45){\ig{5}}\put(25,30){\mor{S}{10}}\put(36,20){=}\egr}
\put(135,10){\bgr{.4ex}{30}{70}\put(0,0){\umor{\eta}{16}}\put(0,24){\coumor{u^*}{16}}\egr}
\cegr

In the third equation of the first row we use equation \rref{sigmaantipode}, in the third equation of the second row we use the equation \rref{rsm}, in the last equation we use \rref{eqrseta}.\eee
 
\begin{cor}\label{koru-1u} $u^*$ is the $*$-inverse of $u$.\end{cor}
\proof
According to \rref{eqdefu}, \rref{eqdefu1} and  \rref{equu}, we have:
 
\cbgr{190}{55}
\put(0,25){\bgr{.4ex}{15}{35}\put(5,0){\rmor{r^*}}\put(5,10){\ig{10}}\put(15,10){\mor{S}{10}}\put(5,20){\comult}\put(0,19){\coumor{u}{6}}\put(0,25){\comult}\put(5,30){\ig{5}}\put(20,10){=}\egr}
\put(23,10){\bgr{.4ex}{30}{60}\put(0,0){\rmor{r^*}}\put(10,10){\ig{10}}\put(0,10){\mor{S^2}{10}}\put(0,20){\comult}\put(5,25){\isibrmor}\put(5,35){\comult}\put(10,40){\comult}\put(15,45){\ig{5}}\put(15,5){\rmor{r^*}}\put(15,15){\mor{S^2}{10}}\put(25,15){\mor{S}{20}}\put(25,35){\twist{-5}{5}}\put(30,25){=}\egr}
\put(58,0){\bgr{.4ex}{30}{60}\put(0,10){\rmor{r^*}}\put(0,20){\mor{S^2}{10}}\put(0,30){\braid}\put(0,40){\twist{5}{5}}\put(5,45){\comult}\put(10,50){\comult}\put(15,55){\ig{5}}\put(10,20){\isibrmor}\put(20,0){\rmor{r^*}}\put(20,10){\mor{S^2}{10}}\put(30,10){\mor{S}{30}}\put(20,30){\ig{10}}\put(10,40){\comult}\put(30,40){\twist{-10}{10}}\put(35,35){=}\egr}
\put(98,5){\bgr{.4ex}{30}{60}\put(0,10){\rmor{r^*}}\put(0,20){\mor{S^2}{10}}\put(0,30){\braid}\put(10,20){\isibrmor}\put(20,0){\rmor{r^*}}\put(20,10){\mor{S^2}{10}}\put(30,10){\mor{S^2}{20}}\put(20,30){\comult}\put(0,40){\comult}\put(5,45){\comult}\put(10,50){\ig{5}}\put(25,35){\twist{-10}{10}}\put(35,30){=}\egr}
\put(140,5){\bgr{.4ex}{30}{60}\put(0,0){\rmor{r^*}}\put(0,10){\ig{10}} \put(20,0){\rmor{r^*}}\put(10,10){\isibrmor}

\put(0,20){\ibraid}
\put(0,30){\comult}\put(5,35){\mor{S}{10}}\put(5,45){\comult}\put(10,50){\ig{5}}\put(15,45){\twist{10}{-10}}\put(20,30){\comult}\put(20,20){\ig{10}}\put(30,10){\mor{S}{20}}\put(34,30){=}\egr}
\put(184,15){\coumor{\va}{30}}
\cegr
The last equation holds according \rref{eqmrs2}.
According to \rref{equsus}, we have:
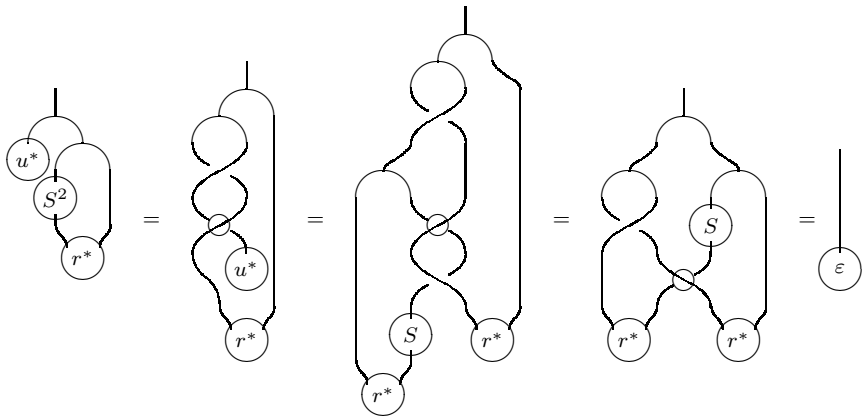
\cbgr{150}{70}
\put(0,25){\bgr{.4ex}{15}{35}\put(5,0){\rmor{r^*}}\put(15,10){\ig{10}}\put(5,10){\mor{S^2}{10}}\put(5,20){\comult}\put(0,19){\coumor{u^*}{6}}\put(0,25){\comult}\put(5,30){\ig{5}}\put(21,10){=}\egr}
\put(30,10){\bgr{.4ex}{30}{60}\put(5,0){\rmor{r^*}}\put(5,10){\twist{-5}{10}}\put(0,20){\sibrmor}\put(0,30){\ibraid}\put(0,40){\comult}\put(5,45){\comult}\put(10,14){\coumor{u^*}{6}}\put(15,10){\ig{35}}\put(10,50){\ig{5}}\put(21,25){=}\egr}
\put(60,0){\bgr{.4ex}{30}{60}\put(0,0){\rmor{r^*}}\put(0,10){\ig{30}}\put(10,10){\mor{S}{10}}\put(0,40){\comult}\put(5,45){\twist{5}{5}}\put(10,50){\ibraid}\put(10,60){\comult}\put(15,65){\comult}\put(20,70){\ig{5}}\put(10,20){\braid}\put(10,30){\sibrmor}\put(20,40){\ig{10}}\put(20,10){\rmor{r^*}}\put(30,20){\ig{40}}\put(30,60){\twist{-5}{5}}\put(36,35){=}\egr}
\put(105,10){\bgr{.4ex}{30}{60}\put(0,0){\rmor{r^*}}\put(20,0){\rmor{r^*}}\put(0,10){\ig{10}}\put(30,10){\ig{20}}\put(10,10){\isibrmor}\put(0,20){\ibraid}\put(0,30){\comult}\put(5,35){\twist{5}{5}}\put(10,40){\comult}\put(15,45){\ig{5}}\put(20,20){\mor{S}{10}}\put(20,30){\comult}\put(25,35){\twist{-5}{5}}\put(36,25){=}\egr}
\put(150,25){\coumor{\va}{25}}
\eeee\cegr

\begin{edl}\label{antip}  Let  $(B,\si,R)$ be a central coquasitriangular Hopf algebra. Then the square of the antipode can be calculated via $u,u^*$ and $\si$:
\bbs S^2= (u\ot B\ot u^*)(B\ot \si\tau^{-1})\Delta^{(2)}.\ees

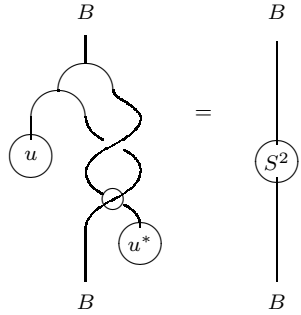
\cbgr{45}{45}\put(0,10){\bgr{.4ex}{10}{20}
\put(0,20){\comult}\put(5,25){\comult}\put(10,10){\ibraid}\put(10,0){\sibrmor}\put(10,-10){\idgr{10}}\put(0,10){\coumor{u}{10}}\put(20,-6){\coumor{u^*}{6}}\put(10,30){\idgr{5}}\put(30,20){=}\put(45,-10){\mor{S^2}{44}}\put(10,-10){\obju{B}}\put(45,-10){\obju{B}}\put(10,35){\objo{B}}\put(45,35){\objo{B}}\put(15,25){\twist{5}{-5}}\egr}  \cegr
\end{edl}

\proof Since $u^*$ is the right inverse of $u$, we have:

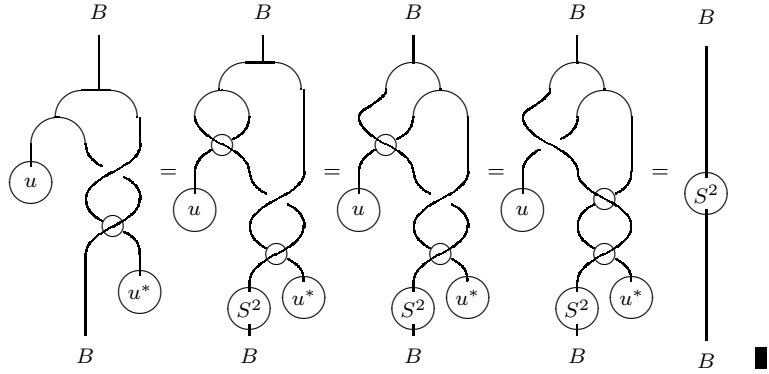
\cbgr{130}{55}\eeee
\put(0,0){\bgr{.4ex}{20}{55}\put(10,0){\idgr{15}}\put(10,15){\sibrmor}
\put(10,25){\ibraid}\put(20,5){\coumor{u^*}{10}}\put(0,25){\coumor{u}{10}}\put(0,35){\comult}\put(5,40){\bcomult{15}}\put(12.5,45){\idgr{10}}\put(20,35){\idgr{5}}
\put(10,0){\obju{B}}\put(12.5,55){\objo{B}}\egr}
\put(25,30){=}
\put(30,0){\bgr{.4ex}{20}{55}\put(10,10){\sibrmor}\put(20,4){\coumor{u^*}{6}}\put(0,20){\coumor{u}{10}}\put(0,30){\isibrmor}\put(0,40){\comult}\put(5,45){\bcomult{15}}\put(12.5,50){\idgr{5}}
\put(20,30){\idgr{15}}\put(10,0){\mor{S^2}{10}}\put(10,20){\ibraid}\put(10,0){\obju{B}}\put(12.5,55){\objo{B}}\egr}
\put(55,30){=}
\put(60,0){\bgr{.4ex}{20}{55}\put(10,20){\ibraid}\put(10,10){\sibrmor}\put(20,4){\coumor{u^*}{6}}\put(0,20){\coumor{u}{10}}\put(0,30){\isibrmor}\put(10,40){\comult}\put(5,45){\comult}\put(10,50){\idgr{5}}
\put(20,30){\idgr{10}}\put(10,0){\mor{S^2}{10}}\put(0,40){\twist{5}{5}}\put(10,0){\obju{B}}\put(10,55){\objo{B}}\egr}
\put(85,30){=}
\put(90,0){\bgr{.4ex}{20}{55}\put(10,20){\isibrmor}\put(10,10){\sibrmor}\put(20,4){\coumor{u^*}{6}}\put(0,20){\coumor{u}{10}}\put(0,30){\braid}\put(10,40){\comult}\put(5,45){\comult}\put(10,50){\idgr{5}}
\put(20,30){\idgr{10}}\put(10,0){\mor{S^2}{10}}\put(0,40){\twist{5}{5}}\put(10,0){\obju{B}}\put(10,55){\objo{B}}\egr}
\put(115,30){=}
\put(125,0){\mor{S^2}{54}}\put(125,0){\obju{B}}\put(125,55){\objo{B}}\cegr

According to  \rref{equsus} and since $u^*$ is the left $*$-inverse of $u$,  we have:

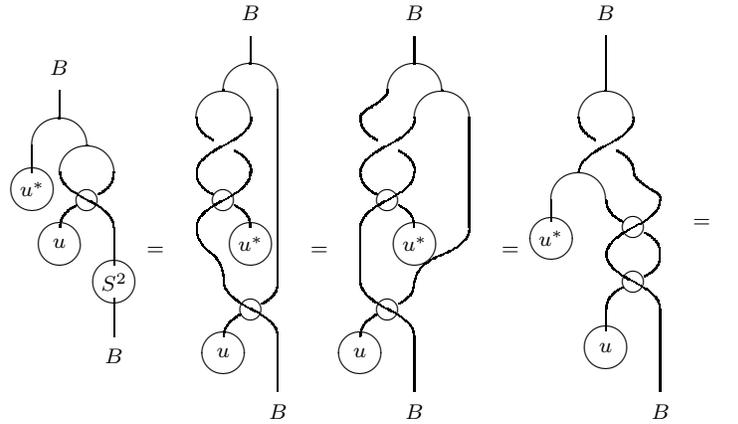
\cbgr{120}{65}
\put(0,10){\bgr{.4ex}{15}{20}\put(0,24){\coumor{u^*}{11}}\put(0,35){\comult}\put(5,40){\ig{5}}\put(5,14){\coumor{u}{6}}\put(5,20){\isibrmor}\put(5,30){\comult}\put(15,0){\mor{S^2}{20}}\put(21,15){=}\put(15,0){\obju{B}}\put(5,45){\objo{B}}\egr}
\put(30,0){\bgr{.4ex}{15}{65}\put(5,4){\coumor{u}{6}}\put(5,10){\isibrmor}\put(15,0){\ig{10}}\put(5,20){\twist{-5}{10}}\put(10,24){\coumor{u^*}{6}}\put(0,30){\sibrmor}\put(0,40){\ibraid}\put(0,50){\comult}\put(5,55){\comult}\put(15,20){\ig{35}}\put(10,60){\ig{5}}\put(21,25){=}\put(15,0){\obju{B}}\put(10,65){\objo{B}}\egr}
\put(60,0){\bgr{.4ex}{15}{65}\put(0,4){\coumor{u}{6}}\put(0,10){\isibrmor}\put(10,0){\ig{10}}\put(0,20){\ig{10}}\put(10,20){\twist{10}{10}}\put(10,24){\coumor{u^*}{6}}\put(0,30){\sibrmor}\put(0,40){\ibraid}\put(0,50){\twist{5}{5}}\put(5,55){\comult}\put(10,60){\ig{5}}\put(10,50){\comult}\put(20,30){\ig{20}}\put(26,25){=}\put(10,0){\obju{B}}\put(10,65){\objo{B}}\egr} 

\put(95,5){\bgr{.4ex}{15}{65}\put(20,-5){\idgr{15}}\put(10,0){\coumor{u}{10}}\put(10,10){\isimor}\put(10,20){\simor}\put(0,30){\comult}\put(5,35){\ibraid}\put(5,45){\comult}\put(10,50){\idgr{10}}\put(0,20){\coumor{u^*}{10}}\put(20,30){\twist{-5}{5}}\put(20,-5){\obju{B}}\put(10,60){\objo{B}}\put(35,15){\ig{30}}\put(26,25){=}\egr}

\cegr

Thus we have proven the following theorem.
\begin{edl} The  antipode $S$ is an isomorphism:
\bbs S^{-1}=(S\ot \id\ot)(u^*\ot u\ot \id)(\id\ot T^{-1})\Delta^{(2)}.\eee\ees
\end{edl}
\noindent{\bf Remarks.} The method used here is based on Majid's proof of the similar fact for quasitriangular Hopf algebras \cite{majid3}.

\section{Constructing central bialgebras from central coalgebras }\label{sec-exam}
The colimit of a monoidal diagram of central coalgebras is a central bialgebra.   Each central bialgebra can be considered as the colimit of the diagram of coalgebras on it. This construction enables us to obtain easily central commutative and coquasitriangular bialgebras. However, central Hopf algebras cannot be obtained by this construction, since the antipode is not a morphism in $\Zz(\Bb(\C))$ for any choice of the bialgebra structure. Equation \rref{sigmaantipode} seems very mysterious for central Hopf algebras other then those reconstructed from monoidal functors.

Recall that a diagram in a category $\C$ is a functor $\theta$ from a small category $\D$ to $\C$, objects of this diagram are the images of objects of $\D$, morphisms are the images morphisms in $\D$. In other words, a diagram $\Aaa$ consists  of objects $A_i$, indexed by objects $i$ in $\D$, morphisms in $\Aaa$ satisfy the axioms for morphisms in a category. A natural morphism $c$ from $\Aaa$ to an object $C$ in $\C$ is, by definition, a family of morphisms $c_i:A_i\lora C$ which commute with morphisms in $\Aaa$. The pairs $(C,c:\Aaa\lora C)$ form a category. The initial object in this category, if it exists, is called colimit of $\Aaa$. Thus, it is an object, say $A$, with morphisms $a_i:A_i\lora A$, called injections, commuting with morphisms in $\Aaa$ and satisfying the following universal property: for any natural morphism $c:\Aaa\lora C$ there exists uniquely a morphism $h:A\lora C$ such that $c=ha$. A category  is said to be cocomplete if any diagram in it possesses a colimit. %The reader is referred to \cite{maclane} for more details.

If $\C$ is a monoidal category, we can define the tensor product of two diagrams $\theta : \D\lora \C$ and $\eta:\E\lora C$ to be the diagram $\ot(\theta\times\eta):\D\times\E\lora \C$. If $\D$ is a (strict) monoidal category and $\theta$ is a monoidal functor then the corresponding diagram is called monoidal. In this case there exist natural isomorphisms $t:I\lora A_1 $ and $t_{i,j}:A_i\ot A_j\lora A_{i\circ j}$, (where $\circ$ denotes the tensor product in $\D$ and $1$ is the unit object in $\D$), satisfying certain coherence-conditions for monoidal functors \cite{maclane}.

We shall call a  category tensor-cocomplete if it is monoidal, cocomplete and the tensor functor  preserves the colimits in both arguments. In this case, for every two diagrams $\Aaa=\{A_i, i\in \D\}$ and $\Bbb=\{B_j, j\in\E\}$ with colimits $(A,a_i)$ and $(B,b_i)$ respectively, their product has the colimit $(A\ot B,a_i\ot b_j)$. We shall prove the following lemma in the Appendix.
\begin{lem}\label{colimit} Let $(\C,\ot)$ be a tensor-cocomplete category.
\begin{itemize}
\item[1)] The category $\Cc(\C)$ of coalgebras in $\C$ is tensor-cocomplete.
\item[2)] The center $\Zz(\C)$ of $\C$ is tensor-cocomplete.
\item[3)] The colimit of a monoidal diagram in $\C$ is an algebra in $\C$.\end{itemize}\end{lem}
Examples of tensor-cocomplete categories are the categories of (co) modules over bialgebras.

Now assume that $(\C,\ot,\tau)$ is a braided tensor-cocomplete category. According to Lemma \ref{colimit}, $\Zz(\Cc(\C))$ is a braided tensor-cocomplete category. If $\Aaa=\{A_i,i\in\D\}$ is a monoidal diagram in $\Zz(\Cc(\C))$ then its colimit $B$ is an algebra in $\Zz(\Cc(\C))$, hence, a central bialgebra in $\C$. In particular, given a central coalgebra $C$, the tensor algebra on it, defined as the colimit of the diagram $\Ccc$, consisting of objects $I,C,C^{\ot 2},\ldots$ and identity-morphisms, is a central bialgebra.

If a monoidal diagram $\Aaa$ contains the braiding $\si$ of $\Zz(\Cc(\C))$ as its morphism then its colimit is a central commutative bialgebra. A weaker condition, assuming that $\Aaa$ is a braided diagram, does not imply the existence of a CQT structure on it colimit.

A bicharacter on a monoidal diagram  $\Aaa=(A_i,i\in\D)$ is a $*$-invertible natural morphism {\it in} $\C:\ r:\Aaa\ot\Aaa\lora I$ subject to the following equations:
\bb\label{bicharacter}\bbar{ll} r_{i,j\circ k}= & r_{i,k}(A_i\ot r_{i,j}\ot A_k)(\Delta_{A_i}\ot A_i\ot A_k),\\
r_{i\circ j,k}= & (r_{i,k}\ot r_{j,k})(A_i\ot \si_{A_j,A_k}\ot A_k)(A_i\ot A_j\ot \Delta_{A_k}).\eear\ee
Let $A$ be the colimit of a monoidal category $\Aaa$ in $\Zz(\Cc(\C))$, then $A$ is central. $A$ is coquasitriangular iff $\Aaa$ possesses a bicharacter, in this case we have:
\bbas r(a_i\ot a_j)=r_{i,j}.\eeas

\begin{lem}\label{bicharecter}
Assume that $r:\Aaa\ot\Aaa\lora I$ is a bicharacter, then $R_{i,j}$ given by
\bb\label{eqRr} R_{i,j}:=(r^{-1}_{i,j}\ot\si_{A_i,A_j}\ot r_{i,j})\Delta^{(2)}_{A_i\ot A_j}:A_i\ot A_j\lora A_j\ot A_i\ee
 is a natural isomorphism in $\Aaa$. \end{lem}
\proof  Indeed, $R_{i,j}$ are coalgebra morphism by their definition, $R_{i,j}$ commute with morphism in $\Aaa$, since $r$ is natural, and $R_{i,j}$ satisfy \rref{braidings} by virtue of \rref{bicharacter}.\eee

For a central bialgebra $B$, let us denote by $\Coalg(B)$ the category of coalgebras $C$ in $\Zz(\Cc(\C))$, together with a coalgebra morphism $b_C:C\lora B$, morphisms in $\Coalg(B)$ are those in $\Zz(\Cc(\C))$ which commute with  $b_C$. Then $\Coalg(B)$ is a monoidal category and $B$ is the colimit of the forgetful functor of this category into $\Zz(\Cc(\C))$. 
\begin{cor} Let $B$ be a central CQT  bialgebra, then  $\Coalg(B)$ is a braided category.\end{cor}

$B$, itself, is an object in $\Coalg(B)$ with $b_B=\id_B$. Thus we have:
\begin{cor} Let $B$ be a central CQT bialgebra, then
\bbas  R_{B,B}=   (r^{-1}\ot\si\ot r)\Delta^{(2)}_{B\ot B}:B\ot B\lora B\ot B\eeas
is a Yang--Baxter operator.\eee\end{cor}
Now let $V$ be a rigid object in $\C$, that is, there exist an object $V^*$ and morphisms $\ev,\db$: $\ev:V^*\ot V\lora I$, $\db:I\lora\ V\ot V^*$ object to the equations:
\bbas V=(V\ot \ev)(\db\ot V),\\
V^*=(\ev\ot V^*)(V^*\ot \db).
\eeas
$V^*\ot V$ is then a coalgebra, the coproduct and counit are $V^*\ot\db\ot V$ and $\ev$ respectively. We can make $V^*\ot V$ into an object of $\Zz(\Cc(\C))$ defining
\bbas \si_{V^*\ot V}(N):=(\tau^{-1}_{V^*,N}\ot N)(V^*\ot\tau_{V,N}).\eeas
Note that $(V^*\ot V,\si_{V^*\ot V})$ is also an object of $\Zz(\C)$.

Let $B$ be a central bialgebra, a rigid $B$-comodule $V$ is called central if the morphism $V^*\ot V\lora B$, induced from the coaction, is in $\Zz(\Cc(\C))$, where $V^*\ot V$ is considered as object of $\Zz(\Cc(\C))$ in the way explained above. Central comodule over a central bialgebra form a category which is braided if the bialgebra is coquasitriangular and rigid (i.e., all its objects are rigid) if the bialgebra is a Hopf algebra.

If a bialgebra $B$ is reconstructed form a monoidal functor $\omega$ from a small category $\V$ into the category $\C_0$ of rigid objects in $\C$ \cite{majid1}, then $B$ is a central bialgebra and, for all $X\in\V,\ \omega(X)$ are central $B$-comodules. Indeed, let us consider a diagram $\Aaa$ in $\Zz(\Cc(\C))$ defined as follows. Objects of $\Aaa$ are the Coends of  diagrams in $\V$ of the following form: $\omega(f):\omega(X)\lora \omega(Y)$, denoted by $\coend_f(X,Y)$, for every pair $ X,Y\in \V$. There only morphism between $\coend_{\id}(X,X)$ and $\coend_f(X,Y)$, and between $\coend_{\id}(Y,Y)$ and $\coend_f(X,Y)$, are induced by the morphism $f:X\lora Y$. Then $B$ is the \hyphenation{co-li-mit} colimit of $\Aaa$. Since $\coend_f(X,Y)$ are objects of $\Zz(\C)$, so is $B$ (cf. \cite{neuchl}).

{\bf Acknowledgment.} The author would like to thank Professor B.~Pareigis for pointing him to this problem and useful advises. He would also like to thank Dr. P.~Schauenburg and M.~Neuchl for interesting discussions. Finally the  author would like to thank the referee for useful comments and suggestions.

This work was supported  in part by the Deutsche Forschungsgemeinschaft and the International Center for Theoretical Physics.

\begin{appendix}\label{appendix}
\section*{Appendix: The proof of Lemma \ref{colimit}}
\subsection*{1. To show that the category $\Cc(\C)$ is cocomplete}
Let $\Aaa=\{A_i,i\in\D\}$ be a diagram  in $\Cc(\C)$ and $A$ be its colimit with the injections $a_i:A_i\lora A$. Let $\Delta_i,\va_i$ denote the coproduct and counit on $A_i$ respectively. The morphism $(a_i\ot a_i)\Delta_i: A_i\lora A\ot A$ is a natural morphism from $\Aaa$ to $A\ot A$, hence induces a morphism $\Delta :A\lora A\ot A$, such that $(a_i\ot a_i)\Delta_i=\Delta a_i$. Analogously induces the natural morphism $\va_i: A_i\lora I$ a morphism $\va:A\lora I$. One has to check:
\begin{itemize}
\item[a)] $(A,\Delta, \va)$ is a coalgebra,
\item[b)] $a_i$ are coalgebra morphisms,
\item[c)] if $c:\Aaa\lora C$ is a natural morphism in $\Cc(\C)$ then there exists a coalgebra morphisms $h:A\lora C$ such that $c=ha$.\end{itemize}

Let us consider the morphism
\bbas j_i=a_i^{\ot 3}\Delta^{(2)}_i:A_i\lora A\ot A\ot A,\eeas
which is a natural morphism from $\Aaa$ to $A\ot A\ot A$. We have
\bbs\bbar{cl}j_i&= (a_i\ot a_i\ot a_i)(\Delta_i\ot A_i)\Delta_i\\
 &=((a_i\ot a_i)\Delta_i\ot a_i)\Delta=(\Delta a_i\ot a_i)\Delta_i\\
 &=(\Delta \ot A)(a_i\ot a_i)\Delta_i= (\Delta\ot A)\Delta a_i.\eear\ees
Analogously one can show $j_i=(A\ot \Delta)\Delta a_i$. From the universal property of $A$ one has $(\Delta\ot A)\Delta=(A\ot \Delta)\Delta$. The assertion for $\va$ can be proven similarly.

b) follows immediately from the definition of $\Delta$ and $\va$.

We now prove c). Let $c:\Aaa\lora C$ be a natural morphism and $k:A\lora C$ be its factor $c=ka$. We have to show $\Delta_C k=(k\ot k)\Delta$. Again we have 
\bbs\bbar{cl} \Delta_Ck a_i&=\Delta_C c_i=(c_i\ot c_i)\Delta_i\\
 & = (k\ot k)(a_i\ot a_i)\Delta_i\\
 &=(k\ot k)\Delta a_i,\eear\ees
whence c) follows.

\subsection*{2. To show that $\Zz(\C)$ is cocomplete}

Let $\Aaa=\{A_i,\si_i, i\in \D\}$ be a diagram in $\Zz(\C)$. Considering it as a diagram in $\C$, we assume that $A$ is its colimit. We have to show:
\begin{itemize}
\item[a)]  $A$ is an object in $\Zz(\C)$,
\item[b)]  $A$ is the colimit of $\{A_i\}$ in $\Zz(\C)$.\end{itemize}
Let $a_i:A_i\lora A$ be the injections. For an object $N$ in $\C$, the morphisms $(N\ot a_i)\si_i(N):A_i\ot N\lora N\ot A_i$ commute with morphisms in $\Aaa$, hence induce morphism $\si_A(N):A\ot N\lora N\ot A$ satisfying:
\bbs (N\ot a_i)\si_i(N)=\si_A(N)(a_i\ot N),\ees
which makes $A$ an object of $\Zz(\C)$.

For showing b) it is enough to show that every natural morphism $b$ from $\{A_i,\si_i\}$ to an object $B$ in $\Zz(\C)$ factors through a morphism $k:A\lora B$, which is in $\Zz(\C)$, that is, for all $N$ in $\C$,
\bbs (N\ot k)\si_A(N)=\si_B(N)(k\ot N).\ees
The morphism $k:A\lora B$ is defined by the universality of $A$, hence satisfies:
\bbs\bbar{cl} \si_B(N)(k\ot N)(a_i\ot N)&= \si_B(N)(k_i\ot N)\\
 & =(N\ot k_i)\si_i(N)\\
 & =(N\ot k)(N\ot a_i)\si_i(N)\\
 &= (N\ot k)\si_A(N)(a_i\ot N).\eear\ees

\subsection*{3. To show that the colimit of a monoidal diagram is an algebra}

Let $\Aaa=\{A_i, i\in \D\}$ be a monoidal diagram and $A$ be its colimit. The naturality of morphisms $t_{i,j}: A_i\ot A_j\lora A_{i\circ j}$ ensures that $a_{i\circ j}t_{i,j}:A_i\ot A_j\lora A$ is a natural morphism from $\Aaa\ot \Aaa$ to $A$, which induces a product on  $ A$. The coherence-conditions for $t_{i,j}$ provides the associativity of this product. The morphism $\eta_A=a_1t:I\lora A$ is the unit of $A$, provided by coherence-conditions for $t$.

\end{appendix}

\end{document}